\documentclass[final,reqno]{siamltex}
\usepackage{latexsym,amsmath,amssymb,amsfonts,mathrsfs}
\usepackage{epsf,graphicx,epsfig,color,cite,cases}
\usepackage{subfigure,graphics,multirow,marginnote,graphicx}
\usepackage{enumerate}
\usepackage{dsfont}
\newcommand{\R}{{\mathbb R}}

\newcommand{\N}{{\mathbb N}}

\newcommand{\Sp}{{\mathbb S}}
\newcommand{\ds}{\displaystyle}
\newcommand{\no}{\nonumber}
\newcommand{\be}{\begin{eqnarray}}
\newcommand{\ben}{\begin{eqnarray*}}
\newcommand{\en}{\end{eqnarray}}
\newcommand{\enn}{\end{eqnarray*}}
\newcommand{\ba}{\backslash}
\newcommand{\pa}{\partial}

\newcommand{\ov}{\overline}

\newcommand{\I}{{\rm Im}}
\newcommand{\Rt}{{\rm Re}}
\newcommand{\g}{\gamma}
\newcommand{\G}{\Gamma}

\newcommand{\vep}{\varepsilon}

\newcommand{\om}{\omega}

\newcommand{\wid}{\widetilde}

\newcommand{\se}{\setminus}

\newcommand{\ol}{\overline}

\newcommand{\half}{\frac{1}{2}}

\newtheorem{remark}[theorem]{Remark}
\newtheorem{algorithm}{Algorithm}[section]

\definecolor{hw}{rgb}{0,0,0}

\begin{document}
\renewcommand{\theequation}{\arabic{section}.\arabic{equation}}

\title{\bf A novel integral equation for scattering by locally rough surfaces
and application to the inverse problem: the Neumann case}
\author{Fenglong Qu \thanks{School of Mathematics and Informational
Science, Yantai University, Yantai, Shandong, 264005, China ({\tt fenglongqu@amss.ac.cn})}
\and Bo Zhang\thanks{NCMIS, LSEC and Academy of Mathematics and Systems Science, Chinese Academy of
Sciences, Beijing, 100190, China and School of Mathematical Sciences, University of Chinese
Academy of Sciences, Beijing 100049, China ({\tt b.zhang@amt.ac.cn})}
\and Haiwen Zhang\thanks{NCMIS and Academy of Mathematics and Systems Science, Chinese Academy of Sciences,
Beijing 100190, China ({\tt zhanghaiwen@amss.ac.cn})}
}

\date{}

\maketitle


\begin{abstract}
This paper is concerned with direct and inverse scattering by a locally perturbed infinite
plane (called a locally rough surface in this paper) on which a Neumann boundary condition is imposed.
A novel integral equation formulation is proposed for the direct scattering problem which is defined on a
bounded curve (consisting of a bounded part of the infinite plane containing the local perturbation and
the lower part of a circle) with two corners and some closed smooth artificial curve.
It is a nontrivial extension of our previous work on direct and inverse scattering by a locally rough
surface from the Dirichlet boundary condition to the Neumann boundary condition
[{\em SIAM J. Appl. Math.}, 73 (2013), pp. 1811-1829].
For the Dirichlet boundary condition, the integral equation obtained is uniquely solvable in the space of bounded continuous functions on the
bounded curve, and it can be solved efficiently by using the Nystr\"om method with a graded mesh.
However, the Neumann condition case leads to an integral equation which is solvable in the space of squarely integrable functions on the bounded curve
rather than in the space of bounded continuous functions, making the integral equation very difficult
to solve numerically. In this paper, we make us of the recursively compressed inverse preconditioning
(RCIP) method developed by Helsing to solve the integral equation which is efficient and capable of
dealing with large wave numbers.
For the inverse problem, it is proved that the locally rough surface is uniquely determined from a
knowledge of the far-field pattern corresponding to incident plane waves.
Further, based on the novel integral equation formulation, a Newton iteration method is developed
to reconstruct the locally rough surface from a knowledge of multiple frequency far-field data.
Numerical examples are also provided to illustrate that the reconstruction algorithm is stable and
accurate even for the case of multiple-scale profiles.

\begin{keywords}
Integral equation, locally rough surface, inverse scattering problem, Neumann boundary condition, RCIP method, Newton iteration method.
\end{keywords}

\begin{AMS}
35R30, 35Q60, 65R20, 65N21, 78A46
\end{AMS}
\end{abstract}

\pagestyle{myheadings}
\thispagestyle{plain}
\markboth{F. Qu, B. Zhang, and H. Zhang}
{A novel method for scattering by rough surfaces}

\section{Introduction}\label{sec1}

This paper is concerned with the problem of scattering of plane acoustic or electromagnetic waves by a
locally perturbed, perfectly reflecting, infinite plane (which is called a locally rough surface).
Such problems arise in many applications such as geophysics, radar, medical imaging, remote sensing
and nondestructive testing (see, e.g., \cite{BaoGaoLi2011,BaoLin2011,BurkardPotthast2010,CZ98,DeSanto1}).

In this paper we are restricted to the two-dimensional case by assuming that the local perturbation
is invariant in the $x_3$ direction. Precisely, let $\G:=\{(x_1,x_2)\;|\;x_2=h(x_1),x_1\in\R\}$ be
the locally rough surface with a smooth function $h\in C^2(\R)$ having a compact support in $\R$.
Denote by $D_+=\{(x_1,x_2)\;|\;x_2>h(x_1),x_1\in\R\}$ the unbounded domain above the
surface $\G$ which is filled with a homogeneous medium. Denote by $k=\omega/c>0$ the wave number
of the wave field in $D_+$, where $\omega$ and $c$ are the wave frequency and speed, respectively.
We assume throughout the paper that the incident field $u^i$ is the plane wave
\ben
u^i(x;d)=\exp({ikd\cdot x}),
\enn
where $d=(\sin\theta,-\cos\theta)\in\Sp^1_-$ is the incident direction, $\theta$ is the angle of
incidence with $-\pi/2<\theta<\pi/2,$ and $\Sp^1_-:=\{x=(x_1,x_2)\;|\;|x|=1,x_2<0\}$ is the lower
part of the unit circle $\Sp^1:=\{x\in\R^2\;|\;|x|=1\}$. Notice that the incident wave is time-harmonic
($e^{-i\om t}$ time dependence), so that the total field $u$ satisfies the Helmholtz equation
\be\label{eq1}
\Delta u+k^2u=0\quad\mbox{in}\;\;D_+.
\en
Here, the total field $u=u^i+u^r+u^s$ satisfies the Neumann boundary condition on the surface $\G$:
\be\label{eq2}
\frac{\pa u}{\pa \nu}=\frac{\pa u^i}{\pa \nu}+\frac{\pa u^r}{\pa \nu}
+\frac{\pa u^s}{\pa \nu}=0\qquad\mbox{on}\;\;\G,
\en
where $\nu$ is the unit normal vector on $\G$ directed into $D_+$, $u^r$ is the reflected wave
of $u^i$ by the infinite plane $x_2=0$:
\ben
u^r(x;d)=\exp(ik[x_1\sin\theta+x_2\cos\theta])
\enn
and $u^s$ is the unknown scattered wave to be determined which is required to satisfy
the Sommerfeld radiation condition
\be\label{rc}
\lim_{r\to\infty}r^{\frac12}\left(\frac{\pa u^s}{\pa r}-iku^s\right)=0,\quad r=|x|,\quad x\in D_+.
\en
This problem models electromagnetic scattering by a locally perturbed, perfectly reflected,
infinite plane in the TM polarization case; it also models acoustic scattering
by a one-dimensional sound-hard, locally rough surface (see Figure \ref{fig4_nr} for the geometric
configuration of the scattering problem). Further, it can be shown that $u^s$ has the following
asymptotic behavior at infinity (see Remark \ref{re3}):
\ben
u^s(x;d)=\frac{e^{ik|x|}}{\sqrt{|x|}}\left(u^\infty(\hat{x};d)+O\Big(\frac{1}{|x|}\Big)\right),
\qquad |x|\to\infty
\enn
uniformly for all observation directions $\hat{x}\in\Sp^1_+$ with
$\Sp^1_+:=\{x=(x_1,x_2)\;|\;|x|=1,x_2>0\}$ the upper part of the unit circle $\Sp^1$.
Here, $u^\infty(\hat{x};d)$ is called the far-field pattern of the scattered field $u^s$,
depending on the incident direction $d\in\Sp^1_-$ and the observation direction $\hat{x}\in\Sp^1_+$.

Direct scattering problems by locally rough surfaces have been studied both numerically and
mathematically. For the Dirichlet case, the well-posedness of the scattering problem was first
proved in \cite{Willers1987} by the integral equation method.
In \cite{BaoLin2011}, the scattering problem is studied by a variational method, based on a
Dirichlet-to-Neumann (DtN) map, and then solved numerically with using the boundary element method.
In \cite{ZZ2013}, a novel integral equation defined on a bounded curve is proposed
for the Dirichlet scattering problem, leading to a fast numerical algorithm for the scattering problem,
even for large wave numbers. However, for the Neumann case, few results are available.
For the special case when the local perturbation is below the infinite plane (which is called the cavity
problem), the well-posedness was established in \cite{ABW2001} via the variational method for the
direct scattering problem with both Neumann and Dirichlet boundary conditions.
A symmetric coupling method of finite element and boundary integral equations wasn developed
in \cite{BaoGaoLi2011}, which can be applied to arbitrarily shaped and filled cavities
with Neumann or Dirichlet boundary conditions.
Recently, in \cite{BHY2018}, the scattering problem by a locally perturbed interface is studied
by a variational method coupled with a boundary integral equation method and numerically solved,
based on the finite element method in a truncated bounded domain coupled with the boundary element method.
It should be mentioned that some studies related to the scattering problem (\ref{eq1})-(\ref{rc}) have
also been conducted extensively. We refer to \cite{Li2012} for cavity scattering problems of
Maxwell's equations, and to \cite{CZ98,CRZ99,CZ99,CM05,CHP06a,CHP06b,CE10,ZC03}
for the non-local perturbation case which is called the rough surface scattering problem.

The inverse scattering problem of reconstructing the rough surfaces has also attracted many
researchers' attention. For example, many numerical algorithms have been developed for
inverse scattering by locally rough surfaces with Dirichlet boundary conditions
(see, e.g., \cite{BaoLin2011,CSS18,CGHIR,DLLY,KressTran2000,ZZ2013,ZZ2017} and the references quoted).
A marching method based on the parabolic integral equation was proposed in \cite{CSS18} for the
reconstruction of a locally rough surface with Dirichlet or Neumann boundary conditions
from phaseless measurements of the single frequency scattering field at grazing angles,
while a Kirsch-Kress decomposition method was given in \cite{LSZ17} to recover a locally rough interface
on which transmission boundary conditions are satisfied.
For numerical recovery of non-locally rough surfaces, we refer to
\cite{BL13,BL14,BaoZhang16,BP2009,BurkardPotthast2010,CSS18,CL05,DeSanto1,DeSanto2,LS14,LZZ18,Spivack}.
For inverse cavity scattering, the reader is referred to \cite{BaoGaoLi2011,Ma05,Li2012}.

In this paper, we extend our previous work in \cite{ZZ2013} on direct and inverse scattering by
a locally rough surface from the Dirichlet condition to the Neumann condition.
Precisely, we propose a novel integral equation formulation for the direct scattering
problem (\ref{eq1})-(\ref{rc}), which is defined on a bounded curve (consisting of a bounded part
of the infinite plane containing the local perturbation and the lower part of a circle) with two
corners and some closed smooth artificial curve.
This extension is nontrivial because, {\color{hw}for the Dirichlet condition case considered in \cite{ZZ2013},
the integral equation obtained in \cite{ZZ2013}
is uniquely solvable in the space of bounded continuous functions on the bounded curve,
and it can be solved efficiently by using the Nystr\"om method with a graded mesh.
However, the Neumann boundary condition leads to an integral equation which is not solvable in the space of bounded continuous functions
on the bounded curve (see, e.g., \cite{bremer2012fast,BOT2009}).}
In this paper, the integral equation formulation obtained is proved to be uniquely solvable
in the space of squarely integrable functions on the bounded curve. But this result then leads to
another difficulty in solving the integral equation numerically since the Nystr\"om method with a
graded mesh used previously in \cite{ZZ2013} does not work anymore.
Instead, we make us of the recursively compressed inverse preconditioning (RCIP) method developed
by Helsing to solve the integral equation (see \cite{Helsing2009,Helsing2012,HelsingKarlsson2013,Helsing2018})
which is fast, accurate and capable of dealing with large wave numbers.
In addition, a further difficulty in our analysis arises in the study of the property of the derived
integral equation formulation in the neighborhood of the two corners of the bounded curve.
To overcome this difficulty, we will use the operator boundedness from \cite{G.A.Chandler1984} in combination
with some detailed energy estimates.

For the associated inverse scattering problem, based on the mixed reciprocity relation established
in this paper, we will prove that the locally rough surface is uniquely determined from a knowledge
of the far-field pattern associate with incident plane waves. Further, we develop a Newton iteration
method to reconstruct the locally rough surface from a knowledge of multiple frequency far-field data.
It should be mentioned that the proposed novel integral equation plays an essential role in solving
the direct scattering problem in each iteration.

This paper is organized as follows.
In Section \ref{sec2}, we first derive a novel integral equation formulation
for the scattering problem (\ref{eq1})-(\ref{rc}) and then prove its unique solvability.
In section \ref{se3}, a fast numerical algorithm, based on the RCIP method, is proposed for solving
the novel integral equation, and a numerical example is carried out to illustrate the performance
of the algorithm. Section \ref{sec3+} is devoted to the proof of uniqueness in the inverse scattering
problem. In Section \ref{se5}, a Newton iteration algorithm with multi-frequency far-field data is
developed to numerically solve the inverse problem, based on the proposed novel integral equation.
Numerical experiments will also be provided to demonstrate that the reconstruction algorithm is stable
and accurate even for the case of multiple-scale profiles. Finally, we will give some conclusion
remarks in Section \ref{se6}.

\section{Solvability of the direct problem via the integral equation method}\label{sec2}
\setcounter{equation}{0}

Let $f=-{\pa u^i}/{\pa\nu}-{\pa u^r}/{\pa\nu}$ on $\G$. It is easy to see that $f$ has a compact support on $\G$.
Then the scattering problem (\ref{eq1})-(\ref{rc}) can be reformulated as the {\bf Neumann problem (NP)}:
Find $u^s\in H^1_{loc}(D_+)$ which satisfies the Helmholtz equation (\ref{eq1}) in $D_+,$
the Sommerfeld radiation condition (\ref{rc}) and the Neumann boundary condition:
\be\label{dbc}
\frac{\pa u^s}{\pa\nu}=f \quad\mbox{on}\quad\G.
\en

\begin{figure}[ht]
\centering
\includegraphics[width=4in]{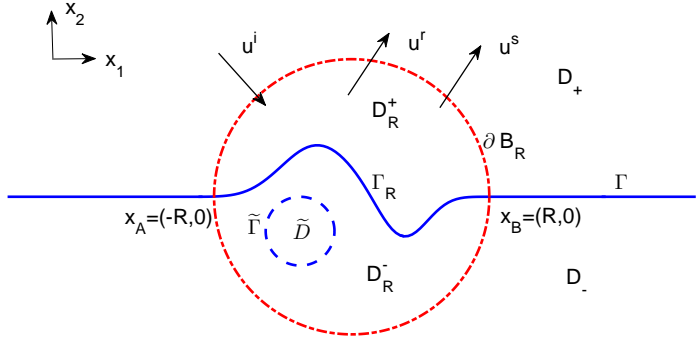}
\vspace{-0.4in}
\caption{The scattering problem from a locally rough surface}\label{fig4_nr}
\end{figure}

We will propose a novel boundary integral equation formulation for the Neumann problem (NP)
which is defined on bounded curves and consequently can be solved with cheap computational cost.
The boundary integral equation formulation will then be used to establish the well-posedness of the
Neumann problem (NP). We first introduce some notations. Assume that $\G_0:=\{(x_1,x_2)\in\G\;|\;x_2=0\}$.
Let $B_R:=\{x=(x_1,x_2)\in\R^2\;|\;|x|<R\}$ be a disk with $R>0$ large enough such that
the local perturbation $\G_p:=\G\ba\G_0=\{(x_1,h(x_1))\;|\;x_1\in\textrm{supp}(h)\}\subset B_R$.
Then $\G_R:=\G\cap B_R$ represents the part of $\G$ containing the local perturbation $\G_p$ of the infinite plane.
Let $x_A:=(-R,0),x_B:=(R,0)$ denote the endpoints of $\G_R$.
For $x=(x_1,x_2)\in\R^2$, let $x^{re}:=(x_1,-x_2)$ be the reflection of $x$ about the $x_1$-axis.
Write $\R^2_\pm:=\{(x_1,x_2)\in\R^2\;|\;x_2\gtrless0\}$,
$D^\pm_R:= B_R\cap D_\pm$ and $\pa B^\pm_R:=\pa B_R\cap D_\pm$,
where $D_-:=\{(x_1,x_2)\;|\;x_2<h(x_1),x_1\in\R\}.$ Further, let $\wid{D}$ be an auxiliary domain in
the interior of $D^-_R$ with a smooth boundary $\wid{\G}$.
See FIG. \ref{fig4_nr} for the geometry of the scattering problem and some notations introduced.

We now have the following uniqueness result for the Neumann problem (NP).

\begin{theorem}\label{thm1-uni}
The Neumann problem (NP) has at most one solution $u^s\in H^1_{loc}(D_+)$.
\end{theorem}

\begin{proof}
Let $u^s$ be the solution of the Neumann problem (NP) with the boundary data $f=0$.
Define the function $\wid{u}^s$ in $\R^2\ba\ov{B}_R$ as follows: $\wid{u}^s(x)=u^s(x)$
for $x\in\R^2_+\ba\ov{B}_{R}$, $\wid{u}^s(x)=u^s(x^{re})$ for $x\in\R^2_-\ba\ov{B}_{R}$.
Then, by the Neumann boundary condition $\pa u^s/\pa\nu=0$ on $\G\se\ov{B}_R$ and the reflection principle
we know that $\wid{u}^s$ satisfies the Helmholtz equation (\ref{eq1}) in $\R^2\ba\ol{B}_R$ and
the Sommerfeld radiation condition (\ref{rc}) uniformly for all directions $\hat{x}=x/|x|$
with $x\in\R^2\ba\ov{B}_{R}$.

Now, let $R_0>R$. By the radiation condition (\ref{rc}) it follows that
\be\label{200}
\int_{\pa B^+_{R_0}}\left\{\left|\frac{\pa u^s}{\pa\nu}\right|^2+k^2|u^s|^2
+2k\textrm{Im}\left(u^s\frac{\pa \ol{u}^s}{\pa\nu}\right)\right\}ds
=\int_{\pa B^+_{R_0}}\left|\frac{\pa u^s}{\pa\nu}-iku^s\right|^2ds
\rightarrow 0\no\\
\en
as $R_0\to+\infty.$ On the other hand, using Green's theorem and the Neumann boundary condition on $\G$ yields
\be\label{201}
\int_{\pa B^+_{R_0}}u^s\frac{\pa \ol{u}^s}{\pa\nu}ds
=2\int_{B_{R_0}\cap D_+}\left(|\nabla u^s|^2-k^2|u^s|^2\right)dx.
\en
Taking the imaginary part of (\ref{201}), we conclude from (\ref{200}) that
\ben
\lim\limits_{R_0\rightarrow+\infty}\int_{\pa B_{R_0}}
\left(\left|\frac{\pa\wid{u}^s}{\pa\nu}\right|^2+k^2|\wid{u}^s|^2\right)ds
=\lim\limits_{R_0\rightarrow+\infty}\int_{\pa B^+_{R_0}}
\left(\left|\frac{\pa u^s}{\pa\nu}\right|^2+k^2|u^s|^2\right)ds=0.
\enn
Thus, by Lemma 2.12 in \cite{ColtonKress2013} we have $\wid{u}^s=0$ in $\R^2\ba\ov{B}_R$,
and so, $u^s=0$ in $\R^2_+\ba\ov{B}_R$. This, together with the unique continuation principle
\cite[Theorem 8.6]{ColtonKress2013}, implies that $u^s=0$ in $D_+$.
The proof is thus complete.
\end{proof}

We will derive a novel integral equation formulation for the Neumann problem (NP) to study
the existence of solutions to the problem (NP) and to develop a fast algorithm to solve the
problem numerically.
To this end, we introduce the layer potentials on the curves $\G_R$, $\pa B^-_R$ and $\wid{\G}$.
For the curve ${I_p}\in\{\G_R,\pa B^-_R,\wid{\G}\}$, define the single-layer potential $\mathcal{S}_{k,I_p}$ by
\be\label{2.2}
(\mathcal{S}_{k,I_p}{\varphi_{I_p}})(x):=\int_{I_p}\Phi_k(x,y){\varphi_{I_p}}(y)ds(y),
\quad x\in\R^2\ba{I_p},
\en
where $\Phi_k(x,y):=(i/4)H^{(1)}_0(k|x-y|),x,y\in\R^2,x\neq y$, is the fundamental solution of
the Helmholtz equation $\Delta w+k^2w=0$ in $\R^2$.
If the wave number $k=0$, $\Phi_0(x,y)=-(2\pi)^{-1}\ln{|x-y|}$ is the
fundamental solution of the Laplace equation in $\R^2$.
By the asymptotic behavior of the fundamental solution $\Phi_k$ (see, e.g., \cite{ColtonKress2013}),
the far-field pattern $S^\infty_{k,I_p}\varphi_{I_p}$ of the single-layer potential
$\mathcal{S}_{k,I_p}\varphi_{I_p}$ can be defined as
\be\label{2.3}
(S^\infty_{k,{I_p}}{\varphi_{I_p}})(\hat{x}):=\frac{e^{i\pi/4}}{\sqrt{8\pi k}}
\int_{I_p}e^{-ik\hat{x}\cdot y}\varphi_{I_p}(y)ds(y), \quad\hat{x}\in\Sp^1.
\en

For ${I_q}\in\{\G_R,\pa B^-_R,\wid{\G}\}$, define the boundary integral operators
$S_{k,{I_p}\to{I_q}}, K'_{k,{I_p}\to{I_q}}$:
\be\label{2.4}
(S_{k,{I_p}\to{I_q}}{\varphi_{I_p}})(x)&:=&\int_{I_p}\Phi_k(x,y){\varphi_{I_p}}(y)ds(y),\quad x\in {I_q},\\
(K'_{k,{I_p}\to{I_q}}{\varphi_{I_p}})(x)&:=&\int_{I_p}\frac{\pa\Phi_k(x,y)}{\pa\nu(x)}{\varphi_{I_p}}(y)ds(y),
\quad x\in {I_q},\label{2.4+}
\en
where $\nu$ is the unit outward normal on $\pa D^-_R$ and $\wid{\G}$.
Further, define the boundary operator $K'^{,re}_{k,{I_p}\to{\pa B^-_R}}$:
\ben
(K'^{,re}_{k,{I_p}\to{\pa B^-_R}}{\varphi_{I_p}})(x)
:=\int_{I_p}\frac{\pa\Phi_k(x^{re},y)}{\pa\nu(x^{re})}{\varphi_{I_p}}(y)ds(y),\quad x\in\pa B^-_R,
\enn
where $\nu$ is the unit outward normal on $\pa B^+_R$ directing into $\R^2\ba\ov{B}_R$.
Some properties of the boundary integral operators $S_{k,{I_p}\to{I_q}}$, $K'_{k,{I_p}\to{I_q}}$ and
$K'^{,re}_{k,{I_p}\to{\pa B^-_R}}$ for ${I_p},{I_q}\in\{\G_R,\pa B^-_R,\wid{\G}\}$
are presented in Remark \ref{re4} below.
For a comprehensive discussion of mapping properties of the layer potentials,
we refer to \cite{ColtonKress2013,MW-2000} and the references therein.

\begin{remark}\label{re4} {\rm
(i) Let ${I_p}\in\{\G_R, \pa B^-_R\}$ and $\varphi_{I_p}\in L^2(I_p)$.
Since any function in $L^2(I_p)$ can be extended by zero to a function in $L^2(\pa D^-_R)$,
it follows from the mapping properties of the single-layer potential over a
closed piecewise-smooth boundary (see \cite[Chapter 6]{MW-2000}) that
$\mathcal{S}_{k,I_p}{\varphi_{I_p}}\in H^{3/2}_{loc}(\R^2)$, satisfies the Helmholtz equation (\ref{eq1})
in $\R^2\ba\ov{I_p}$ and the jump relations on the curve $I_p$:
\be\label{eq43}
(\mathcal{S}_{k,{I_p}}{\varphi_{I_p}})_\pm|_{I_p}&=&S_{k,{I_p}\to{I_p}}{\varphi_{I_p}},\\ \label{eq44}
\left(\frac{\pa}{\pa\nu}\mathcal{S}_{k,{I_p}}{\varphi_{I_p}}\right)_\pm\Big|_{I_p}
&=&\mp\half\varphi_{I_p} +K'_{k,{I_p}\to{I_p}}{\varphi_{I_p}},
\en
where
\ben
(\mathcal{S}_{k,{I_p}}{\varphi_{I_p}})_\pm(x)&:=&\lim\limits_{h\to0^+}(\mathcal{S}_{k,{I_p}}
{\varphi_{I_p}})(x\pm h\nu(x)),\;\;x\in I_p,\\
\lim\limits_{h\to 0^+}\left(\frac{\pa}{\pa\nu}\mathcal{S}_{k,{I_p}}{\varphi_{I_p}}\right)_\pm(x)
&:=&\left(\nu\cdot\nabla\mathcal{S}_{k,{I_p}}{\varphi_{I_p}}\right)(x\pm h\nu(x)),\;\;x\in I_p.
\enn
For ${I_q}\in\{\G_R,\pa B^-_R,\wid{\G}\}$ it is easy to see that
$S_{k,{I_p}\to{I_q}}:L^2(I_p)\to H^1(I_q)$, $K'_{k,{I_p}\to{I_q}}:L^2(I_p)\to L^2(I_q)$
and $K'^{,re}_{k,{I_p}\to{\pa B^-_R}}:L^2(I_p)\to L^2({\pa B^-_R})$ are bounded operators.

(ii) Let $\varphi_{\wid{\G}}\in H^{-1/2}(\wid{\G})$. Then, from \cite[Chapter 6]{MW-2000} it follows that
$\mathcal{S}_{k,\wid{\G}}\varphi_{\wid{\G}}\in H^1_{loc}(\R^2)$ and satisfies the Helmholtz equation (\ref{eq1})
in $\R^2\ba\ov{\wid{\G}}$ as well as the jump relations (\ref{eq43}) and (\ref{eq44}) on the curve $I_p=\wid{\G}$.
Further, since $\wid{\G}$ is a $C^2$-smooth curve, by \cite{ColtonKress2013} it is known that
$S_{k,\wid{\G}\to\wid{\G}}$ and $K'_{k,\wid{\G}\to\wid{\G}}$ are compact operators from
$H^{-1/2}(\wid{\G})$ to $H^{-1/2}(\wid{\G}))$.
Moreover, it can be seen that $K'_{k,\wid{\G}\to\G_R}:H^{-1/2}(\wid{\G})\to L^2(\G_R)$
and $K'_{k,\wid{\G}\to\pa B^-_R},K'^{,re}_{k,\wid{\G}\to\pa B^-_R}:H^{-1/2}(\wid{\G})\to L^2(\pa B^-_R))$
are bounded operators since the kernels of those operators are smooth functions.
}
\end{remark}

Let $u^s$ be the solution to the Neumann problem (NP).
We now extend $u^s(x)$ into $\R^2_-\ba\ov{B}_R$ by reflection,
which we denote by $u^s(x)$ again, such that $u^s(x)=u^s(x^{re})$ for $x\in\R^2_-\ba\ov{B}_R$.
Then, from the reflection principle and the Neumann boundary condition
$\pa u^s/\pa\nu=-({\pa u^i}/{\pa\nu}+{\pa u^r}/{\pa\nu})=0$ on $\G\se\ov{B}_R$
it follows that $u^s\in H^1_{loc}(\R^2\ba\ov{D^-_R})$ and
satisfies the Helmholtz equation (\ref{eq1}) in $\R^2\ba\ov{D_R^-}$.
Following the idea in \cite{bremer2012fast}, we seek the solution $u^s$ in the form
\be\no
u^s(x)&=&\int_{\G_R}\Phi_k(x,y)\varphi_1(y)ds(y)+\int_{\pa B^-_R}\Phi_k(x,y)\varphi_2(y)ds(y)
+\int_{\wid{\G}}\Phi_k(x,y)\varphi_3(y)ds(y),\\ \label{2.6}
&&\qquad\hspace{6cm} x\in\R^2\ba(\pa D_R^-\cup\wid{\G})
\en
with $\Phi=(\varphi_1,\varphi_2,\varphi_3)^T\in X:=L^2(\G_R)\times L^2(\pa B^-_R)\times H^{-\frac{1}{2}}(\wid{\G})$
with the norm $\|\Phi\|_X:=\|\varphi_1\|_{L^2(\G_R)}+\|\varphi_2\|_{L^2(\pa B^-_R)}
+\|\varphi_3\|_{H^{-1/2}(\wid{\G})}$.
We will prove that, by choosing appropriate auxiliary curve $\wid{\G}$
the Neumann problem (NP) can be reduced to an equivalent boundary integral equation
which is uniquely solvable (see Theorem \ref{thm1} and Remark \ref{re6} below).

Let $\psi^{re}$ be a continuous mapping from $\pa D_R^-$ to $\pa D_R^+$ such that
\ben
\psi^{re}(x)=\left\{\begin{array}{ll}
   x, &x\in\Gamma_R,\\
   x^{re}, &x\in\pa B_R^-\cup\{x_A,x_B\}.
   \end{array}\right.
\enn
Then, by the boundary condition (\ref{dbc}) and the reflection principle, $u^s(x)$ satisfies
\be\label{eq27}
\frac{\pa u^s(x)}{\pa\nu(x)}=-\left(\frac{\pa u^i(x)}{\pa\nu(x)}+\frac{\pa u^r(x)}{\pa\nu(x)}\right)
      &&\quad\mbox{for}\;x\in\G_R,\\ \label{eq40}
\frac{\pa u^s(x)}{\pa\nu(x)}-\frac{\pa u^s(\psi^{re}(x))}{\pa\nu(x^{re})}=0
      &&\quad\mbox{for}\;x\in\pa B_R^-.
\en
We now impose the impedance boundary condition for $u^s$ on the auxiliary curve $\wid{\G}$:
\be\label{eq28}
\frac{\pa u^s(x)}{\pa\nu(x)}+i\rho u^s(x)=0\quad\mbox{for}\;x\in\wid{\G},
\en
where $\rho$ is a positive constant. Then the jump relations (see Remark \ref{re4}) and boundary conditions
(\ref{eq27})-(\ref{eq28}) lead to the boundary integral equation
\be\label{2.7}
(I+A)\Phi=G,
\en
where $I$ is the identity operator and $A,G$ are given by
\be\label{eq10}
A:=(A_{ij})_{3\times 3}
\en
with $A_{11}=-2K'_{k,\G_R\to\G_R}$, $A_{12}=-2K'_{k,\pa B^-_R\to\G_R}$, $A_{13}=-2K'_{k,\wid{\G}\to\G_R}$,
$A_{21}=-2(K'_{k,\G_R\to\pa B^-_R}-K^{\prime,re}_{k,\G_R\to\pa B^-_R})$,
$A_{22}=-2(K'_{k,\pa B^-_R\to\pa B^-_R}-K^{\prime,re}_{k,\pa B^-_R\to\pa B^-_R})$,
$A_{23}=-2(K'_{k,\wid{\G}\to\pa B^-_R}-K^{\prime,re}_{k,\wid{\G}\to\pa B^-_R})$,
$A_{31}=-2(K'_{k,\G_R\to\wid{\G}}+i\rho S_{k,\G_R\to\wid{\G}})$,
$A_{32}=-2(K'_{k,\pa B^-_R\to\wid{\G}}+i\rho S_{k,\pa B^-_R\to\wid{\G}})$,
$A_{33}=-2(K'_{k,\wid{\G}\to\wid{\G}}+i\rho S_{k,\wid{\G}\to\wid{\G}})$
and
\be\label{eq11}
G:=\left(\begin{array}{c}
    2\left({\pa u^i}/{\pa\nu}+{\pa u^r}/{\pa\nu}\right)\big|_{\G_R} \\
     0\\
     0 \\
  \end{array}\right).
\en
Obviously, $G\in X$. Further, it is seen from Remark \ref{re4} that $A$ is a bounded linear operator on $X$.

Conversely, we have the following result.

\begin{lemma}\label{le1}
Let $u^s$ be given by (\ref{2.6}) with $\Phi=(\varphi_1,\varphi_2,\varphi_3)^T\in X$
which is the solution of the boundary integral equation (\ref{2.7}) with $A$ and $G$ defined in (\ref{eq10})
and (\ref{eq11}), respectively. Then $u^s\in H^1_{loc}(\R^2\ba\ov{D_R^-})$ and solves the Neumann problem (NP).
\end{lemma}

\begin{proof}
By the fact that $\Phi=(\varphi_1,\varphi_2,\varphi_3)^T\in X$ and Rermark \ref{re4},
it is easy to see that $u^s$ defined by (\ref{2.6}) satisfies the Helmholtz equation (\ref{eq1})
in $\R^2\ba\ov{D_R^-}$ and the Sommerfeld radiation condition (\ref{rc}) and is in $H^1_{loc}(\R^2\ba\ov{D_R^-})$.
Further, by the integral equation $(I+A)\Phi=G$ and the jump relations of the single-layer potentials
(see Rermark \ref{re4}), it follows that
${\pa u^s(x)}/{\pa\nu(x)}=-{\pa u^i(x)}/{\pa\nu(x)}-{\pa u^r(x)}/{\pa\nu(x)}$ for $x\in\G_R$ and
${\pa u^s(x)}/{\pa\nu(x)}={\pa u^s(x^{re})}/{\pa\nu(x^{re})}$ for $x\in\pa B_R^-$.

Let $\wid{u}^s(x)=u^s(x^{re})$ for $x\in\R^2\ba\ov{B}_R$. Then $\wid{u}^s$ satisfies the Helmholtz
equation (\ref{eq1}) in $\R^2\ba\ov{B}_R$, the Sommerfeld radiation condition (\ref{rc})
and the condition ${\pa\wid{u}^s(x)}/{\pa\nu(x)}={\pa u^s(x)}/{\pa\nu(x)}$ for $x\in\pa B_R$.
Further, it follows from the uniqueness of the exterior Neumann problem
(see, e.g., \cite[Chapter 3]{ColtonKress2013}) that $\wid{u}^s=u^s$ in $\R^2\ba\ov{B}_R$,
which implies that $u^s(x)=u^s(x^{re})$ in $\R^2\ba\ov{B}_R$.
In particular, we obtain that ${\pa u^s}/{\pa\nu}=0$ on $\G\ba\G_R$.
Therefore, we have ${\pa u^s}/{\pa\nu}=-({\pa u^i}/{\pa\nu}+{\pa u^r}/{\pa\nu})$ on $\G$.
The proof is thus complete.
\end{proof}

The following theorem gives the unique solvability of the integral equation $(I+A)\Phi=G$.

\begin{theorem}\label{thm1}
Assume that $k^2$ is not a Dirichlet eigenvalue of $-\Delta$ in $\wid{D}$.
Let $A$ and $G$ be given by (\ref{eq10}) and (\ref{eq11}), respectively.
Then the integral equation $(I+A)\Phi=G$ has a unique solution $\Phi=(\varphi_1,\varphi_2,\varphi_3)^T\in X$
with the estimate
\be\label{2.8}
||\varphi_1||_{L^2(\G_R)}+||\varphi_2||_{L^2(\pa B_R^-)}+||\varphi_3||_{H^{-1/2}(\wid{\G})}
\leq C\left\|\frac{\pa u^i}{\pa\nu}+\frac{\pa u^r}{\pa\nu}\right\|_{L^2(\G)}
\en
\end{theorem}

\begin{proof}
The proof is broken down into the following steps.

{\bf Step I.} We show that $I+A$ is a Fredholm operator of index zero.

{\bf Step I.1.} Define
\ben
&&M_0:=\left(M_{0ij}\right)_{3\times 3}\nonumber\\
&&\;=\left(\begin{array}{lll}
    -2K'_{0,\G_R\to\G_R}& -2K'_{0,\pa B^-_R\to\G_R} & 0\\
     -2\left(K'_{0,\G_R\to\pa B^-_R}-K^{\prime,re}_{0,\G_R\to\pa B^-_R}\right)&
     -2\left(K'_{0,\pa B^-_R\to\pa B^-_R}-K^{\prime,re}_{0,\pa B^-_R\to\pa B^-_R}\right)& 0 \\
     0 & 0 & 0
  \end{array}\right).
\enn
Then it follows from Remark \ref{re4} that $M_0$ is bounded in $X$. In addition,
the kernels of $A_{ij}-M_{0ij}$, $i,j=1,2$, are bounded with the upper bound
\ben
C\left[\max\left(\ln\frac{1}{|x-y|},1\right)
+\max\left(\ln\frac{1}{|\psi^{re}(x)-y|},1\right)\right]
\enn
for $x,y\in\pa D^-_R$ with $x\neq y$ and thus weakly singular since ${|x-y|}\leq C{|\psi^{re}(x)-y|}$
for $x,y\in\pa D^-_R$.
Therefore, $A_{ij}-M_{0ij}$, $i,j=1,2$, are compact.
Further, the operators $A_{13}, A_{23}, A_{31}, A_{32}$ are compact since the kernels of those operators
are smooth functions, and, by Remark \ref{re4}, the operator $A_{33}$ is also compact.
Consequently, $A-M_0$ is compact in $X$.

For $z\in\R^2$ and $r\in\R$ define $B_r(z):=\{x\in\R^2\;|\;|x-z|<r\}$.
Obviously, we can choose a fixed $r_0>0$ such that for $z=x_A,x_B$, $B_r(z)\cap(\G\ba\G_0)=\emptyset$
and $\ov{B_r(x_A)}\cap\ov{B_r(x_B)}=\emptyset$ for $r\leq r_0$.
In the remaining part of the proof, we will choose $r\leq r_0$.
Now, choose a cut-off function $\psi_{r,z}\in C^\infty_0(\R^2)$ satisfying that $0\leq\psi_{r,z}\leq 1$,
$\psi_{r,z}(x)=1$ in the region $0\leq|x-z|\leq r/2$ and $\psi_{r,z}(x)=0$ in the region $|x-z|\ge r$.
Define $M_{0,r}:X\to X$ by $M_{0,r}:=\left((M_{0,r})_{ij}\right)_{3\times3}$ with
\ben
(M_{0,r})_{ij}\varphi_j:=\left\{
\begin{array}{ll}
\psi_{r,x_A}(M_0)_{ij}(\psi_{r,x_A}\varphi_j)
+\psi_{r,x_B}(M_0)_{ij}(\psi_{r,x_B}\varphi_j),& i,j=1,2\\
      0,& \textrm{otherwise}
\end{array}\right.
\enn
for any $\Phi=(\varphi_1,\varphi_2,\varphi_3)^T\in X$.
Since the kernel of $M_0-M_{0,r}$ vanishes in a neighborhood of $(x_A,x_A)$ and $(x_B,x_B)$,
$M_0-M_{0,r}$ is compact. Then the operator $A-M_{0,r}=(A-M_0)+(M_0-M_{0,r})$ is compact from $X$ into $X$.

{\bf Step I.2.} We show that $||M_{0,r}||_{X\to X}<1$ for a sufficiently small constant $r>0$.

By the definition of $M_{0,r}$, we have that for any $\Phi=(\varphi_1,\varphi_2,\varphi_3)^T\in X$,
\ben
\|M_{0,r}\Phi\|_X&=&\left\|\sum_{j=1}^2\left(\psi_{r,x_A}(M_0)_{1j}(\psi_{r,x_A}\varphi_j)
+\psi_{r,x_B}(M_0)_{1j}(\psi_{r,x_B}\varphi_j)\right)\right\|_{L^2(\G_R)}\\
&&+\left\|\sum_{j=1}^2\left(\psi_{r,x_A}(M_0)_{2j}(\psi_{r,x_A}\varphi_j)
+\psi_{r,x_B}(M_0)_{2j}(\psi_{r,x_B}\varphi_j)\right)\right\|_{L^2(\pa B^-_R)}\\
&\leq&\sum_{j=1}^2\|\left(\psi_{r,x_A}(M_0)_{1j}(\psi_{r,x_A}\varphi_j)
+\psi_{r,x_B}(M_0)_{1j}(\psi_{r,x_B}\varphi_j)\right)\|_{L^2(\Gamma_R)}\\
&&
+\sum_{j=1}^2\|\left(\psi_{r,x_A}(M_0)_{2j}(\psi_{r,x_A}\varphi_j)
+\psi_{r,x_B}(M_0)_{2j}(\psi_{r,x_B}\varphi_j)\right)\|_{L^2(\pa B^-_R)}.
\enn
Recalling that $\ov{B_r(x_A)}\cap\ov{B_r(x_B)}=\emptyset$ and $\textrm{supp}(\psi_{r,z})\subset\ov{B_r(z)}$
for $z=x_A,x_B$, we find that
\ben
&&\|M_{0,r}\Phi\|_X\le\sum_{j=1}^2\left(\|\psi_{r,x_A}(M_0)_{1j}(\psi_{r,x_A}\varphi_j)\|^2_{L^2(\G_R)}
+\|\psi_{r,x_B}(M_0)_{1j}(\psi_{r,x_B}\varphi_j)\|^2_{L^2(\G_R)}\right)^{\half}\\
&&\qquad\quad+\sum_{j=1}^2\left(\|\psi_{r,x_A}(M_0)_{2j}(\psi_{r,x_A}\varphi_j)\|^2_{L^2(\pa B^-_R)}
+\|\psi_{r,x_B}(M_0)_{2j}(\psi_{r,x_B}\varphi_j)\|^2_{L^2(\pa B^-_R)}\right)^{\half}.
\enn
We will estimate each term of the right-hand side of the above inequality.
In doing so, the following inequality plays an essential role:
\be\label{eq15}
|\nu(x)\cdot(x-y)|\leq C|x-y|^2
\en
for $x,y\in\G_R$ or for $x\in\pa B_R,\;y\in\pa B_R^-$ (see \cite[Section 3.5]{ColtonKress2013}).

We first estimate the norm of $\psi_{r,z}(M_0)_{ij}(\psi_{r,z}\varphi_j)$, $i,j=1,2$
with $z=x_A$ in the following four steps.

(i) For $\psi_{r,x_A}(M_0)_{11}(\psi_{r,x_A}\varphi_1)$ we have
\ben
\psi_{r,x_A}(M_0)_{11}(\psi_{r,x_A}\varphi_1)&=&-2\psi_{r,x_A}K'_{0,\G_R\to\G_R}(\psi_{r,x_A}\varphi_1)\\
&=&-2\psi_{r,x_A}(x)\int_{\G_R\cap B_r(x_A)}\frac{\pa\Phi_0(x,y)}{\pa\nu(x)}\psi_{r,x_A}(y)\varphi_1(y)ds(y).
\enn
It then follows from (\ref{eq15}) that
\be\label{eq21}
&&\|\psi_{r,x_A}(M_0)_{11}(\psi_{r,x_A}\varphi_1)\|^2_{L^2(\G_R)}\no\\
&&\qquad\qquad\le\|(M_0)_{11}\psi_{r,x_A}\varphi_1\|^2_{L^2(\G_R\cap B_r(x_A))}\no\\
&&\qquad\qquad=\int_{\G_R\cap B_r(x_A)}\left(\int_{\G_R\cap B_r(x_A)}
2\frac{\pa\Phi_0(x,y)}{\pa\nu(x)}\psi_{r,x_A}(y)\varphi_1(y)ds(y)\right)^2ds(x)\no\\
&&\qquad\qquad\le C\int_{\G_R\cap B_r(x_A)}\|\varphi_1\|^2_{L^2(\G_R\cap B_r(x_A))}ds(x)\no\\
&&\qquad\qquad\le Cr\|\varphi_1\|^2_{L^2(\G_R\cap B_r(x_A))}.
\en

(ii) We estimate the norm of $\psi_{r,x_A}(M_0)_{12}(\psi_{r,x_A}\varphi_2)$.
Since $x_A=(-R,0)$, the curve $\G_R\cap B_r(x_A)$ can be parameterized as $x(t)=(-R+t,0)$, $t\in[0,r]$.
Further, it is easily seen that $\pa B^-_R\cap \pa B_r(x_A)=(-\sqrt{R^2-\tau^2_r},-\tau_r)$
with $\tau_r:={r\sqrt{4 R^2-r^2}}/({2R})$. Thus $\pa B^-_R\cap B_r(x_A)$ can be parameterized
as $y(\tau)=(y_1(\tau),y_2(\tau)):=(-\sqrt{R^2-\tau^2},-\tau)$, $\tau\in[0,{r\sqrt{4 R^2-r^2}}/({2R})]$.
Then $y_1(0)=-R, y'_1(0)=0$. Choose $\delta\in(0,1)$. Then there exists a $\tau(\delta)>0$ such that
$|y_1(\tau)-(-R)|\leq\delta\tau$ and $|y'_1(\tau)|\leq\delta$ for all $\tau\in[0,\tau(\delta)]$.
Thus we have that for $t\in[0,r],\;\tau\in[0,\tau(\delta)]$,
\be\label{eq25}
|y'(\tau)|&=&\sqrt{|y'_1(\tau)|^2+|y'_2(\tau)|^2}\leq\sqrt{1+\delta^2}\leq 1+\delta,\\ \label{eq16}
|x(t)-y(\tau)|^2&=&(-R+t-y_1(\tau))^2+\tau^2\no\\
&=&t^2-2t(y_1(\tau)+R)+\tau^2\geq t^2-\delta t^2-\frac{1}{\delta}(y_1(\tau)+R)^2+\tau^2\no\\
&\geq&(1-\delta)t^2-\frac{1}{\delta}(\delta\tau)^2+\tau^2=(1-\delta)(t^2+\tau^2).
\en
Further, from (\ref{eq16}) it follows that
\be\label{eq26}
\frac{|y_2(\tau)|}{|x(t)-y(\tau)|}\leq\frac{1}{1-\delta}\frac{\tau}{t^2+\tau^2}.
\en
Now, choose $r>0$ small enough so that ${r\sqrt{4R^-r^2}}/({2R})\le\min[\tau(\delta),1]$ and
$r\le\min[r_0,1]$. Then, applying the parameterizations of $x(t)$ and $y(\tau)$, and using
the inequalities (\ref{eq25}) and (\ref{eq26}), we obtain that
\ben
&&\|\psi_{r,x_A}(M_0)_{12}(\psi_{r,x_A}\varphi_2)\|^2_{L^2(\G_R)}\\
&&\qquad\le\int_{\G_R\cap B_r(x_A)}\left(\int_{\pa B^-_R\cap B_r(x_A)}
2\frac{\pa\Phi_0(x,y)}{\pa\nu(x)}\psi_{r,x_A}(y)\varphi_2(y)ds(y)\right)^2ds(x)\\
&&\qquad\le\int_0^r\left(2\int_0^{\frac{\sqrt{4R^2-r^2}}{2R}r}
\left|\frac{\pa\Phi_0(x,y)}{\pa\nu(x(t))}\right|\cdot
|\varphi_2(y(\tau))|\cdot|y'(\tau)|d\tau\right)^2|x'(t)|dt\\
&&\qquad\le \int_0^r\left(\frac{1}{\pi}\int_0^{\frac{\sqrt{4R^2-r^2}}{2R}r}\frac{|y_2(\tau)|}{|x(t)-y(\tau)|^2}
|\varphi(y(\tau))|\cdot|y'(\tau)|d\tau\right)^2|x'(t)|dt\\
&&\qquad\le \int_0^r\left(\frac{1+\delta}{1-\delta}\int_0^{\frac{\sqrt{4R^2-r^2}}{2R}r}
\frac{\tau}{\pi(t^2+\tau^2)}|\varphi(y(\tau))|d\tau\right)^2dt.
\enn
Let
\ben
\wid{\varphi}_2(\tau)=\left\{
\begin{array}{ll}
\ds \varphi_2(y(\tau)) &\quad \textrm{if}\quad 0\leq\tau\le{r\sqrt{4R^2-r^2}}/({2R}),\\
\ds 0 &\quad\textrm{if}\quad{r\sqrt{4R^2-r^2}}/({2R})<\tau\le 1.
\end{array}\right.
\enn
Then
\be\label{eq17}
\|\psi_{r,x_A}(M_0)_{12}(\psi_{r,x_A}\varphi_2)\|_{L^2(\G_R)}\leq\frac{1+\delta}{1-\delta}
\left(\int_0^1\left(\int_0^1\frac{\tau}{\pi(t^2+\tau^2)}|\wid{\varphi_2}(\tau)|d\tau\right)^2dt\right)^\half.\no\\
\en
For any $v\in L^2(0,1)$, define the operator
\ben
(Rv)(t):=\int^1_0\frac{t}{\pi(t^2+\tau^2)}v(\tau)d\tau.
\enn
From \cite[Lemma 1]{G.A.Chandler1984} it is known that $R$ is a linear bounded operator
on $L^2(0,1)$ with the norm $\|R\|_{L^2(0,1)\to L^2(0,1)}\leq\sin({\pi}/{4})=1/\sqrt{2}$.
Further, it is easy to see that
\ben
\left(\int_0^1\left(\int_0^1\frac{\tau}{\pi(t^2+\tau^2)}|\wid{\varphi_2}(\tau)|d\tau\right)^2dt\right)^\half
=\|R^*|\wid{\varphi}_2|\|_{L^2(0,1)}.
\enn
It then follows by (\ref{eq17}) and the fact that $\|R^*\|_{L^2(0,1)\to L^2(0,1)}=\|R\|_{L^2(0,1)\to L^2(0,1)}$
that
\be\label{eq18}
\|\psi_{r,x_A}(M_0)_{12}(\psi_{r,x_A}\varphi_2)\|_{L^2(\G_R)}
\leq\frac{1+\delta}{1-\delta}\|R^*|\wid{\varphi}_2|\|_{L^2(0,1)}
\leq\frac{1+\delta}{1-\delta}\frac{1}{\sqrt{2}}\|\wid{\varphi}_2\|_{L^2(0,1)}.\no\\
\en
Now, since $|y'(\tau)|=\sqrt{(y'_1(\tau))^2+(y'_2(\tau))^2}\ge|y'_2(\tau)|= 1$ for $\tau\in[0,\tau(\delta)]$,
we have
\be\label{eq19}
\|\varphi_2\|_{L^2(\pa B^-_R\cap B_r(x_A))}=\left(\int_0^{\frac{\sqrt{4R^2-r^2}}{2R}r}
\left(\varphi_2(y(\tau))\right)^2|y'(\tau)|d\tau\right)^\half\ge\|\wid{\varphi}_2\|_{L^2(0,1)}.\quad\quad
\en
It thus follows from (\ref{eq18}) and (\ref{eq19}) that
\be\label{eq22}
\|\psi_{r,x_A}(M_0)_{12}(\psi_{r,x_A}\varphi_2)\|_{L^2(\G_R)}\leq
\frac{1+\delta}{1-\delta}\frac{1}{\sqrt{2}}\|\varphi_2\|_{L^2(\pa B^-_R\cap B_r(x_A))}.
\en

(iii) For $\psi_{r,x_A}(M_0)_{21}(\psi_{r,x_A}\varphi_1)$, we have
\be\label{eq20}
&&\psi_{r,x_A}(M_0)_{21}(\psi_{r,x_A}\varphi_1)\no\\
&&\quad=-2\psi_{r,x_A}(K'_{0,\G_R\to\pa B^-_R}-K^{\prime,re}_{0,\G_R\to\pa B^-_R})(\psi_{r,x_A}\varphi_1)\no\\
&&\quad=-2\psi_{r,x_A}(x)\int_{\G_R}\left(\frac{\pa\Phi_0(x,y)}{\pa\nu(x)}
-\frac{\pa\Phi_0(x^{re},y)}{\pa\nu(x)}\right)\psi_{r,x_A}(y)\varphi_1(y)ds(y).
\en
It is easy to verify that for $y\in\G_R\cap B_r(x_A)$ and $x\in\pa B^-_R\cap B_r(x_A)$,
\ben
\frac{\pa\ln|x-y|}{\pa\nu(x)}=\frac{\nu(x)\cdot(x_1-y_1,x_2)}{|x-y|^2}.
\enn
Then, by the definition of $x^{re}$ and $\nu(x^{re})$ we have that for $y\in\G_R\cap B_r(x_A)$
and $x\in\pa B^-_R\cap B_r(x_A)$,
\ben
\frac{\pa\ln|x-y|}{\pa\nu(x)}-\frac{\pa\ln|x^{re}-y|}{\pa\nu(x)}=0.
\enn
Combining this with (\ref{eq20}) implies that
\be\label{eq23}
\|\psi_{r,x_A}(M_0)_{21}(\psi_{r,x_A}\varphi_1)\|_{L^2(\pa B^-_R)}=0.
\en

(iv) With the inequality (\ref{eq15}), the estimate of $\psi_{r,x_A}(M_0)_{22}(\psi_{r,x_A}\varphi_2)$
is similar to that of $\psi_{r,x_A}(M_0)_{11}(\psi_{r,x_A}\varphi_1)$, and so we have
\be\label{eq24}
\|\psi_{r,x_A}(M_0)_{22}(\psi_{r,x_A}\varphi_2)\|^2_{L^2(\G_R)}
\leq C r\|\varphi_2\|^2_{L^2(\pa B^-_R\cap B_r(x_A))}.
\en

The norm of $\psi_{r,z}(M_0)_{ij}(\psi_{r,z}\varphi_j)$, $i,j=1,2$ with $z=x_B$ can be estimated similarly,
and we have similar estimates as (\ref{eq21}), (\ref{eq22}), (\ref{eq23}) and (\ref{eq24}).

We are now ready to estimate the norm of $M_{0,r}\Phi$. For arbitrarily fixed $\delta\in(0,1)$,
we obtain by combining the estimates obtained in (i)-(iv) that there exists a $r>0$ small enough
such that
\ben
\|M_{0,r}\Phi\|_X^2&\le& Cr\left(\|\varphi_1\|^2_{L^2(\G_R\cap B_r(x_A))}
   +\|\varphi_1\|^2_{L^2(\G_R\cap B_r(x_B))}\right)\\
&&+\left[Cr+\frac12\left(\frac{1+\delta}{1-\delta}\right)^2\right]
\left(\|\varphi_2\|^2_{L^2(\pa B^-_R\cap B_r(x_A))}+\|\varphi_2\|^2_{L^2(\pa B^-_R\cap B_r(x_B))}\right)\\
&\le& Cr\|\varphi_1\|^2_{L^2(\G_R)}
 +\left[Cr+\frac12\left(\frac{1+\delta}{1-\delta}\right)^2\right]\|\varphi_2\|^2_{L^2(\pa B^-_R)}\\
&\le&\left[Cr+\frac12\left(\frac{1+\delta}{1-\delta}\right)^2\right]\|\Phi\|^2_X.
\enn
Fix $0<\delta<(\sqrt{2}-1)/(\sqrt{2}+1)$ to obtain that $({1+\delta})/({1-\delta})<\sqrt{2}$, so we have
$\|M_{0,r}\|_{X\to X}<1$ for $r>0$ small enough.

Now, since $\|M_{0,r}\|_{X\to X}<1$ for $r>0$ small enough, $I+M_{0,r}$ has a bounded inverse in $X$.
Thus, and since $I+A=I+M_{0,r}+(A-M_0)+(M_0-M_{0,r})$, we know by Steps I.1 and I.2 that
$I+A$ is a Fredholm operator of index zero. The proof of Step I is thus complete.

{\bf Step II.} To prove that $I+A$ is injective and thus invertible in $X$.

Let $(I+A)\Phi=0$ for $\Phi=(\varphi_1,\varphi_2,\varphi_3)^T\in X$.
Then it follows from Lemma \ref{le1} that the scattered field $u^s$,
given by (\ref{2.6}), satisfies the boundary condition ${\pa u^s}/{\pa \nu}=0$ on $\G$.
By Theorem \ref{thm1-uni} and the unique continuation principle we have $u^s=0$ in $\R^2\ba\ol{D_R^-}$.
Then, by the jump relations of the layer potentials (see Remark \ref{re4}) one has
\be\label{2.11}
u^s_-=u^s_+=0\quad\mbox{on}\quad\pa D_R^-.
\en
The integral equation $(I+A)\Phi=0$ implies the boundary condition ${\pa u^s_+}/{\pa\nu}+i\rho u^s_+=0$
on $\wid{\G}$. Thus $u^s$ satisfies the mixed boundary value problem
\be\label{2.12}
\left\{\begin{array}{ll}
       \ds\Delta u^s+k^2u^s=0& \mbox{in}\;\;D^-_R\ba\ov{\wid{D}},\\
       \ds u^s_-=0& \mbox{on}\;\;\pa D_R^-,\\
       \ds \frac{\pa u^s_+}{\pa\nu}+i\rho u^s_+=0 & \mbox{on}\;\;\wid{\G}.\\
\end{array}\right.
\en
From (\ref{2.12}) and Green's formula it follows that
\ben
\int_{D^-_R\ba\ov{\wid{D}}}\left(|\nabla u^s|^2-k^2|u^s|\right)ds=i\rho\int_{\wid{\G}}|u^s_+|^2ds.
\enn
Taking the imaginary part of the above equation gives that $u^s_+=0$ on $\wid{\G}$.
This, together with the last equation in (\ref{2.12}), implies that ${\pa u^s_+}/{\pa\nu}=0$ on $\wid{\G}$.
Holmgren's uniqueness theorem then yields that $u^s=0$ in $D^-_R\ba\ov{\wid{D}}$.

By the jump relations of the layer potentials on $\pa D_R^-$, one has
\ben
-\varphi_1&=&\frac{\pa u^s_+}{\pa\nu}-\frac{\pa u^s_-}{\pa\nu}=0 \quad\mbox{on}\quad \G_R,\\
-\varphi_2&=&\frac{\pa u^s_+}{\pa\nu}-\frac{\pa u^s_-}{\pa\nu}=0 \quad\mbox{on}\quad \pa B_R^-.
\enn
On the other hand, the jump relations of the layer potentials on $\wid{\G}$ give
\ben
u^s_-=u^s_+=0\quad \mbox{on}\quad \wid{\G}.
\enn
We then obtain that $u^s$ satisfies the impedance problem in $\wid{D}$:
\ben
\left\{\begin{array}{ll}
\ds \Delta u^s+k^2u^s=0 & \mbox{in}\;\;\wid{D},\\
\ds u^s_-=0 & \mbox{on}\;\;\wid{\G}.
\end{array}\right.
\enn
Since $k$ is not a Dirichlet eigenvalue of $-\Delta$ in $\wid{D}$, we have $u^s=0$ in $\wid{D}$.
The jump relations of layer potentials on $\wid{\G}$ again imply that
\ben
-\varphi_3=\frac{\pa u^s_+}{\pa\nu}-\frac{\pa u^s_-}{\pa\nu}=0\quad\mbox{on}\quad\wid{\G}.
\enn
Therefore, we have $\Phi=(\varphi_1,\varphi_2,\varphi_3)^T=0$, yielding that $I+A$ is injective on $X$.

By Step I and the Fredholm alternative it follows that $I+A$ has a bounded inverse on $X$, and
the estimate (\ref{2.8}) is thus obtained. This completes the proof of the theorem.
\end{proof}

\begin{remark}\label{re6}{\rm
In Theorem \ref{thm1}, it is assumed that $k^2$ is not a Dirichlet eigenvalue of $-\Delta$ in $\wid{D}$.
It should be noted that in numerical computation this assumption can be easily fulfilled.
In fact, $\wid{D}$ can be chosen to be a ball of radius $r$.
Then, by \cite{ColtonKress2013} (see the proof of Corollary 5.3 of \cite{ColtonKress2013})
we know that $k^2$ is not a Dirichlet eigenvalue of $-\Delta$ in $\wid{D}$
if and only if $kr$ is not the zeros of the Bessel function $J_n$ of order $n$ for any integer $n$.
In particular, if the radius $r$ is chosen so that $kr$ is smaller than the smallest zero of
the Bessel function $J_0$, then $k^2$ is not a Dirichlet eigenvalue of $-\Delta$ in $\wid{D}$.
}
\end{remark}

\begin{remark}\label{re3a} {\rm
It should be remarked that the circle $\pa B_R$ can be replaced by any closed smooth curve
which is symmetric about the plane $x_2=0$.
However, it is computationally simple and convenient to use the circle $\pa B_R$.
}
\end{remark}

We now conclude the following result on the well-posedness of the problem (NP) by
combining Lemma \ref{le1}, Theorems \ref{thm1-uni} and \ref{thm1} and the mapping properties of
the single- and double-layer potentials in $X$.

\begin{theorem}\label{thm1-wp}
The problem (NP) has a unique solution $u^s\in H^1_{loc}(D_+)$. Further, for any $R>0$ we have
\be\label{eq13_nr}
||u^s||_{H^1(D^+_R)}\leq C\left\|\frac{\pa u^i}{\pa\nu}+\frac{\pa u^r}{\pa\nu}\right\|_{L^2(\G)}
\en
\end{theorem}

\begin{remark}\label{re3}{\rm
By the asymptotic behavior of the fundamental solution $\Phi_k$ (see, e.g., \cite{ColtonKress2013})
and the form (\ref{2.6}) of the scattered field $u^s$, it is easy to know that
$u^s$ has the following asymptotic behavior at infinity
\ben
u^s(x;d)=\frac{e^{ik|x|}}{\sqrt{|x|}}\left(u^\infty(\hat{x};d)+O\Big(\frac{1}{|x|}\Big)\right),
\qquad |x|\to\infty
\enn
uniformly for all observation directions $\hat{x}\in S_+$. Here, the far-field pattern
$u^\infty(\hat{x};d)$ is given by
\be\label{eq14_nr}
u^\infty(\hat{x})=S^\infty\Phi:=(S^\infty_{k,\G_R}{\varphi_1})(\hat{x})
+(S^\infty_{k,{\pa B^-_R}}{\varphi_2})(\hat{x})+(S^\infty_{k,\wid{\G}}{\varphi_3})(\hat{x}),
\en
which is an analytic function on the unit circle $\Sp^1$,
where $\Phi=(\varphi_1,\varphi_2,\varphi_3)^T\in X$ is the unique solution of the integral
equation $(I+A)\Phi=G$ with $A$ and $G$ given by (\ref{eq10}) and (\ref{eq11}), respectively.
}
\end{remark}

\begin{remark}\label{re2-9} {\rm
We may seek the solution $u^s$ in the form
\be\label{eq32}
u^s(x)=\int_{\G_R}\Phi_k(x,y)\varphi_1(y)ds(y)+\int_{\pa B^-_R}\Phi_k(x,y)\varphi_2(y)ds(y),
\;\; x\in\R^2\ba\pa D_R^-\qquad
\en
with ${\Phi}=(\varphi_1,\varphi_2)^T\in L^2(\Gamma_R)\times L^2(\pa B^-_R)$.
Then, by using the boundary conditions (\ref{eq27}) and (\ref{eq40}) and the jump relations of layer potentials,
we obtain the boundary integral equation
\be\label{eq31}
(I+B){\Phi}={G},
\en
where ${B}:=({B}_{ij})_{2\times 2}$ with ${B}_{ij}=A_{ij}$, $i,j=1,2$ and $A_{ij}$ given in (\ref{eq10})
and ${G}=(2({\pa u^i}/{\pa\nu}+{\pa u^r}/{\pa\nu})|_{\G_R},0)^T$.

By a similar argument as those used in the proof of Theorem \ref{thm1}, we can prove that
if $k^2$ is not a Dirichlet eigenvalue of $-\Delta$ in $D^-_R$ then $I+B$ is boundedly invertible
on $L^2(\G_R)\times L^2(\pa B^-_R)$ and that if ${\Phi}=(\varphi_1,\varphi_2)^T$
is the unique solution to the integral equation (\ref{eq31}) then
the scattering field $u^s$ in the form (\ref{eq32}) solves the Neumann problem (NP).
Again, by the asymptotic behavior of the fundamental solution $\Phi_k$,
the far-field pattern of the scattered field $u^s$ in the form (\ref{eq32}) is given by
\be\label{eq39}
u^\infty(\hat{x})=(S^\infty_{k,\G_R}{\varphi_1})(\hat{x})+(S^\infty_{k,\pa B^-_R}{\varphi_2})(\hat{x}).
\en
It is seen that the integral equation (\ref{eq31}) does not involve the auxiliary curve $\wid{\G}$ and
is much simpler than the integral equation (\ref{2.7}). In the numerical experiments
for inverse problem we will make use of the integral equation (\ref{eq31})
and the far-field pattern (\ref{eq39}) to solve the solution of the direct scattering problem in
each iteration (see Section \ref{se5} for details).
}
\end{remark}

\section{Numerical solution of the direct scattering problem}\label{se3}

\subsection{The numerical algorithm}\label{se3-1}

In this subsection, we develop a numerical algorithm to solve the direct scattering problem,
employing the boundary integral equation (\ref{2.7}).
Note that the kernels of the integral operators in the integral equation (\ref{2.7}) are hyper-singular
near the corners $x_A$ and $x_B$ of the curve $\pa D^-_R$, leading to difficulties in numerical computation.
To deal with this issue, we make use of the RCIP method to solve the integral equation (\ref{2.7}).
The RCIP method was first proposed by J. Helsing to solve planar boundary value problems in the presence
of various boundary singularities and has been successfully applied to many kinds of static problems
in nonsmooth domains including the Helmholtz equation
(see \cite{Helsing2009,Helsing2012,HelsingKarlsson2013,Helsing2018}).
In what follows, we will present the steps of our numerical algorithm
based on the idea in \cite{HelsingKarlsson2013}. For more details and applications of the RCIP method,
the reader is referred to \cite{Helsing2009,Helsing2012,HelsingKarlsson2013} and the references therein.

For the integral equation (\ref{2.7}), we rewrite $A$ and $\Phi$ as follows
\ben
A=\left(\begin{array}{cc}
    \hat{A}_{11} & \hat{A}_{12} \\
    \hat{A}_{21} & \hat{A}_{22} \\
  \end{array}\right),\;
\Phi=\left(\begin{array}{c}
    \hat{\varphi}_{1}  \\
    \hat{\varphi}_{2}  \\
  \end{array}\right)
\enn
with
\ben
&&\hat{A}_{11}:=\left(\begin{array}{cc}
    A_{11} & A_{12} \\
    A_{21} & A_{22} \\
  \end{array}\right),\;
\hat{A}_{12}:=\left(\begin{array}{c}
    A_{13}  \\
    A_{23} \\
  \end{array}\right),\;
\hat{A}_{21}:=\left(\begin{array}{cc}
    A_{31} & A_{32} \\
  \end{array}\right),\;
\hat{A}_{22}:=A_{33}\\
&&\hat{\varphi}_{1}:=\left(\begin{array}{c}
    \varphi_{1}  \\
    \varphi_{2} \\
  \end{array}\right),\;
\hat{\varphi}_2:=\varphi_3.
\enn
It can be seen that $\hat{A}_{11}$ is the integral operator $B$ on $\pa D^-_R$ defined in (\ref{eq31})
in Remark \ref{re2-9}.
Denote by $\hat{A}_{11}(x,y)$ the kernel of $\hat{A}_{11}$ and
split the kernel $\hat{A}_{11}(x,y)$ into two parts
\be\label{eq42}
\hat{A}_{11}(x,y)=\hat{A}^\star_{11}(x,y)+\hat{A}^\circ_{11}(x,y)
\en
such that $\hat{A}^\star_{11}(x,y)$ is zero except when both $x$ and $y$ are close enough
to the same corner of $\pa D^-_R$. Thus $\hat{A}^\circ_{11}(x,y)$
is zero when $x,y$ are close enough to the same corner of $\pa D^-_R$.
Following the kernel split (\ref{eq42}), we have the corresponding operator split
\ben
\hat{A}_{11}=\hat{A}^\star_{11}+\hat{A}^\circ_{11},
\enn
where $\hat{A}^\circ_{11}$ is a compact operator.
Arguing similarly as in the proof of the fact that $||M_{0,r}||_{X\to X}<1$ for a sufficiently
small constant $r>0$ in Step I.2 of the proof of Theorem \ref{thm1}, we can prove that
$\|\hat{A}^\star_{11}\|<1$ if, in the split of $\hat{A}_{11}$ above, $\hat{A}^\star_{11}(x,y)$ is
zero except when both $x$ and $y$ are close enough to the same corner of $\pa D^-_R$, and so $(I+\hat{A}^\star_{11})^{-1}$ exists. By the variable substitution
\ben
\Phi=\left(\begin{array}{cc}
    (I+\hat{A}^\star_{11})^{-1} &0 \\
   0 & I \\
  \end{array}\right)
\wid{\Phi},\qquad
\wid{\Phi}:=\left(\begin{array}{c}
    \widetilde{\hat{\varphi}}_{1}  \\
    \widetilde{\hat{\varphi}}_{2}  \\
  \end{array}\right),
\enn
the integral equation (\ref{2.7}) is rewritten as a right preconditioned integral equation
\be\label{eq33}
\wid{\Phi}(x)+\left(\begin{array}{cc}
    \hat{A}^\circ_{11} & \hat{A}_{12} \\
    \hat{A}_{21} & \hat{A}_{22} \\
  \end{array}\right)
\left(\begin{array}{cc}
    (I+\hat{A}^\star_{11})^{-1} &0 \\
   0 & I \\
  \end{array}\right)
\wid{\Phi}(x)=G(x).
\en

\begin{figure}[htbp]
\subfigure{\begin{minipage}[b]{0.5\textwidth}
\centering
\includegraphics[width=1\textwidth]{{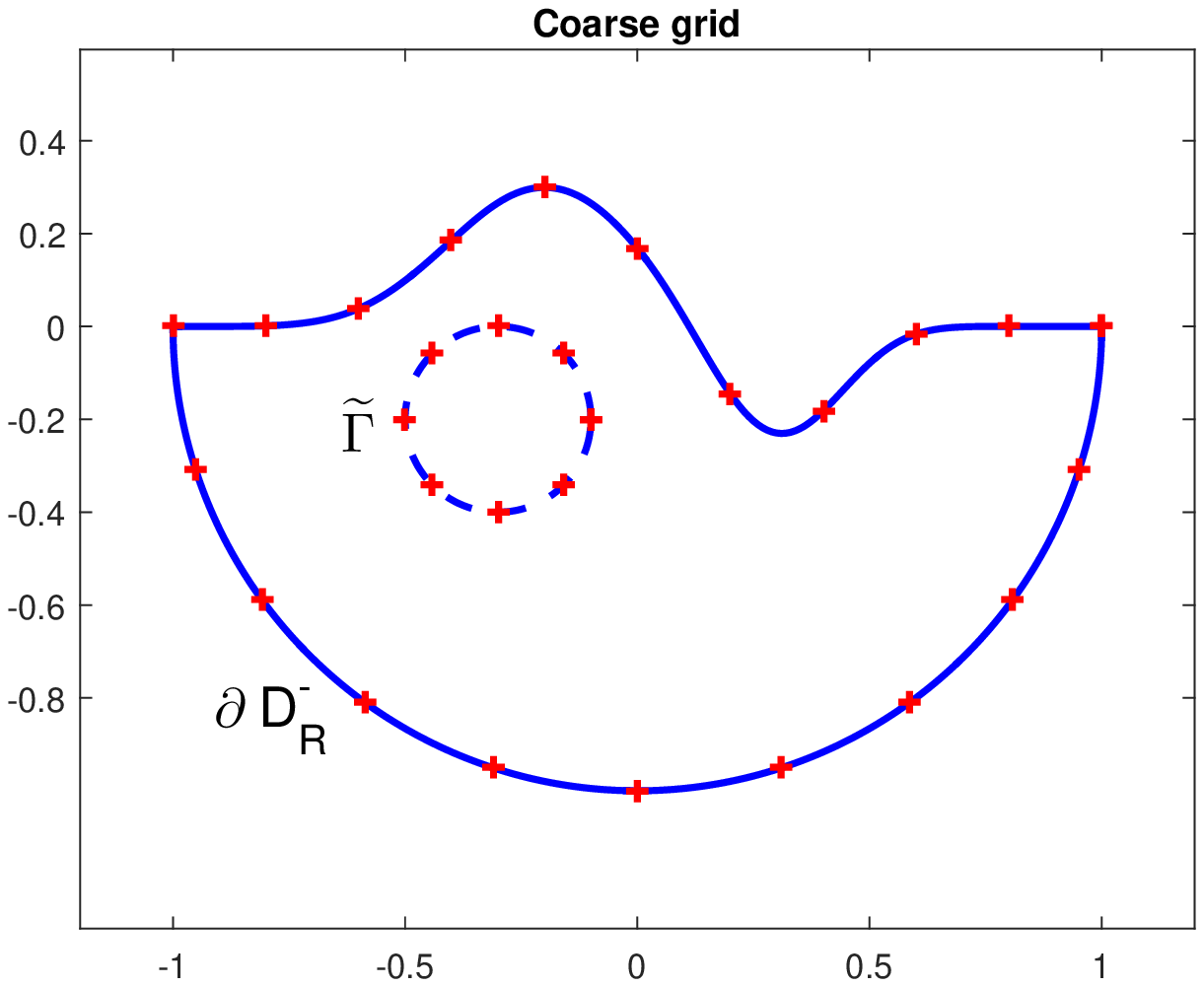}}
\end{minipage}}%
\subfigure{\begin{minipage}[b]{0.5\textwidth}
\centering
\includegraphics[width=1\textwidth]{{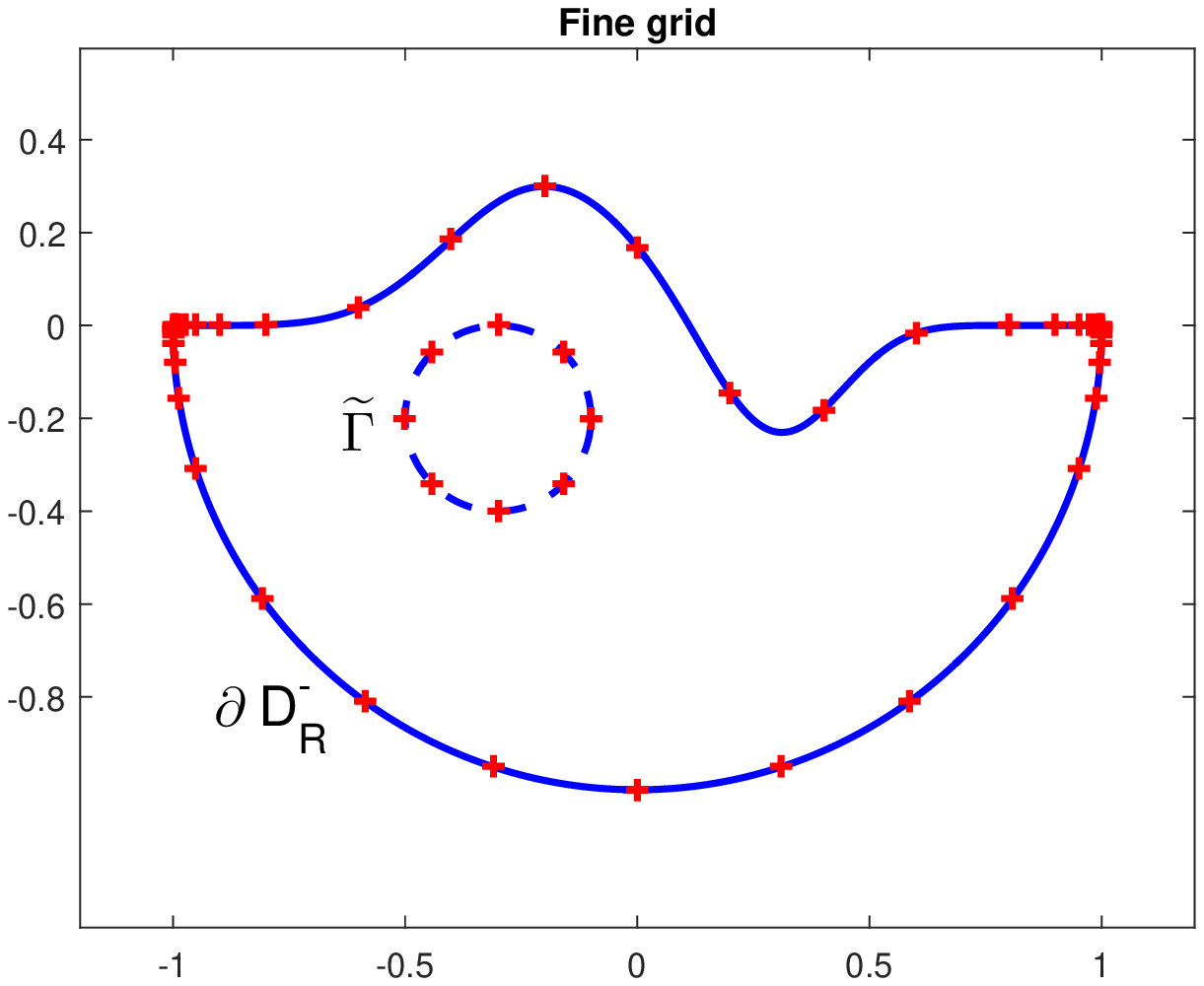}}
\end{minipage}}
\caption{Geometry of the coarse and fine meshes
}\label{fig6}
\end{figure}

To discretize the integral equation (\ref{eq33}) we need two different meshes: a coarse mesh and a fine mesh.
For the coarse mesh, we divide each smooth curve $I_p\in\{\G_R,\pa B^-_R,\wid{\G}\}$ into $n_{pan}$ panels
which have approximately equal length.
The fine mesh is constructed around the corner $x_z$ from the coarse one by halving the two panels around
the corner $x_z$ each time and repeatedly $n_{sub}$ times, $z=A,B$.
Thus, for the curves $\pa D^-_R$ and $\wid{\G}$, there are totally $3n_{pan}$ panels on the coarse mesh
and $3n_{pan}+4 n_{sub}$ panels on the fine mesh.
For the geometry of the coarse and fine meshes, see FIG. \ref{fig6}.

With these two meshes, the integral equation (\ref{eq33}) can be discretized by using a Nystr\"{o}m scheme
based on the composite $16$-point Gauss-Legendre quadrature (see \cite[Section 2]{Helsing2009}).
The quantities $\wid{\Phi},\hat{A}^\circ_{11},\hat{A}_{12},\hat{A}_{21},\hat{A}_{22}$
are simply discretized on the coarse mesh. Only the operator $(I+\hat{A}^\star_{11})^{-1}$ needs
to be discretized on the fine mesh to get an accurate approximation.
Then the discretized version of the integral equation (\ref{eq33}) can be obtained as follows
\be\label{eq34}
\bf\rm
\left[\mathbf{I}_{coa}+\left(
\begin{array}{cc}
  \rm\mathbf{\hat{A}}^\circ_{11,coa} &\rm\mathbf{\hat{A}}_{12,coa} \\
  \rm\mathbf{\hat{A}}_{21,coa} &\rm \mathbf{\hat{A}}_{22,coa} \\
\end{array}\right)
\left(\begin{array}{cc}
  \bf\rm \mathbf{R} &\bf\rm \mathbf{0} \\
  \bf\rm \mathbf{0} &\bf\rm \mathbf{I}_{coa} \\
\end{array}\right)\right]\mathbf{\widetilde{\Phi}}_{coa}=\mathbf{G}_{coa},
\en
where $\rm\mathbf{\wid{\Phi}}_{coa},\mathbf{\hat{A}}^\circ_{11,coa},\mathbf{\hat{A}}_{12,coa}$,
$\mathbf{\hat{A}}_{21,coa},\mathbf{\hat{A}}_{22,coa}$ are the discretized matrixes of $\wid{\Phi}$, $\hat{A}^\circ_{11}$, $\hat{A}_{12}$, $\hat{A}_{21}$, $\hat{A}_{22}$, respectively,
and $\mathbf{R}$ is called the compressed weighted inverse matrix and given by
\ben\bf\rm
\mathbf{R}=\mathbf{P}^T_W(\mathbf{I}_{fin}+\mathbf{\hat{A}}^\star_{11,fin})^{-1}\mathbf{P}.
\enn
Here, the subscript "coa" indicates the grid on the coarse mesh and the subscript "fin" indicates
the grid on the fine mesh. The prolongation matrix $\mathbf{P}$ performs a polynomial interpolation
from the coarse grid to the fine grid and $\rm\mathbf{P}^T_W$ is the transpose of a weighted
prolongation matrix (see \cite[Section 6]{Helsing2012}).
Note that $\mathbf{R}$ depends on the fine mesh but can be easily computed by fast and stable
recursion (see Sections 7, 8, 12 and 13 in \cite{Helsing2012}).

By solving the linear equation (\ref{eq34}), we get the unknown $\rm\mathbf{\widetilde{\Phi}}_{coa}$.
Finally, the far-field pattern defined in (\ref{eq14_nr}) can be approximated as follows:
\ben
u^\infty\backsimeq\rm{\mathbf{S}^\infty_{coa}}
\left(\begin{array}{cc}
   \bf\rm \mathbf{R} &\bf\rm \mathbf{0} \\
   \bf\rm \mathbf{0} &\bf\rm \mathbf{I}_{coa} \\
\end{array}\right)\mathbf{\widetilde{\Phi}}_{coa},
\enn
where ${\rm\mathbf{S}^\infty_{coa}}$ is the discretized matrix of $S^\infty$ on the coarse mesh of also
using the composite 16-point Gauss-Legendre quadrature.

\subsection{Numerical experiments}

In this section, a numerical experiment is carried out to illustrate the performance of
the algorithm proposed in Section \ref{se3-1}. We consider the locally rough surface with
\ben
h(x_1)=\left\{\begin{array}{ll}
\ds \exp\left[16/(25x_1^2-16)\right]\sin(4\pi x_1), &|x_1|<4/5,\\
\ds 0, &|x_1|\geq4/5.
\end{array}\right.
\enn
Generally speaking, there is no exact solution of the scattering problem (NP) when the incident field
is given by an incident plane wave. However, in some special cases, we can find a solution
of the scattering problem (NP) which has an explicit expression. In this example,
we consider the incident field to be the point source $u^i(x;y):=\Phi_k(x,y)$ with $y=(-0.1,0.1)$.
Then the reflected field can be characterized by $u^r(x;y)=\Phi_k(x,y^{re})$.
Suppose $u^s(x;y)$ is the solution to the scattering problem (NP) with the boundary data
$f=-{\pa u^i(x;y)}/{\pa\nu(x)}-{\pa u^r(x;y)}/{\pa\nu(x)}$.
Since $u^i(x;y)+u^r(x;y)$ satisfies the Helmholtz equation (\ref{eq1}) in $D_+$, the Neumann boundary
condition ${\pa u^i(x;y)}/{\pa\nu(x)}+{\pa u^r(x;y)}/{\pa\nu(x)}=0$ on $\G\ba\G_R$ and
the radiation condition (\ref{rc}), then it is easy to see that the scattered field
$u^s(x;y)=-\Phi_k(x,y)-\Phi_k(x,y^{re})$. Thus, it follows from the asymptotic behavior of $\Phi_k$
that the far-field pattern $u^\infty(\hat{x};y)$ of the scattered field $u^s(x;y)$ is given by
\ben
u^\infty(\hat{x};y)=-\frac{e^{i\pi/4}}{\sqrt{8\pi k}}
\left(e^{-ik\hat{x}\cdot y}+e^{-ik\hat{x}\cdot y^{re}}\right),
\enn
which has an explicit expression and can be accurately computed.

We now compute the approximated value $u^\infty_{app}(\hat{x};y)$ of the far-field pattern
$u^\infty(\hat{x};y)$, using the RCIP method proposed in Section \ref{se3-1} with
the following parameters.

1) The wave number $k=100$ and the observation direction $\hat{x}=(\cos(2\pi/3),\sin(2\pi/3))$.

2) The auxiliary curve $\wid{\G}$ is taken to be a disk of radius $r=0.1$ and centered at $[0,-0.5]$.
By Remark \ref{re6}, it is valid to use the integral equation (\ref{2.7}) to solve the Neumann problem (NP).

3) The parameter $\rho$ in the boundary condition (\ref{eq28}) is chosen as $\rho=1$.

4) For the coarse mesh, we choose the number of panels on each smooth component to be $n_{pan}=5N$
with $N=1,2,3,\ldots, 24$, while, for the fine mesh, we choose the number of subdivisions to
be $n_{sub}=30$.

\begin{figure}[htbp]
\centering
\includegraphics[width=3in]{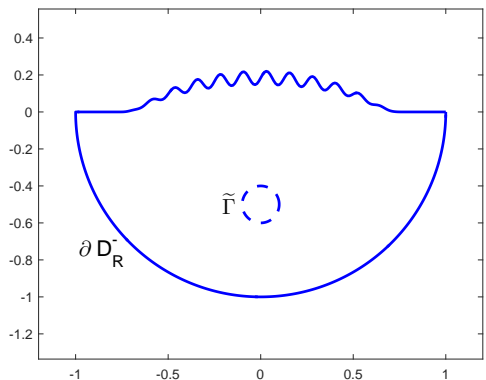}
\caption{Geometry of the problem
}\label{fig7}
\end{figure}

For the geometry of the problem, see FIG. \ref{fig7}.

\begin{figure}[htbp]
\centering
\includegraphics[width=3in]{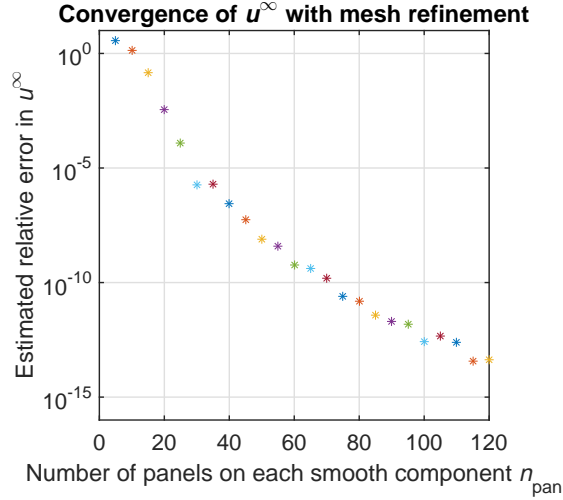}
\caption{The relative error ${|u^\infty(\hat{x};y)-u^\infty_{app}(\hat{x};y)|}/{|u^\infty(\hat{x};y)|}$
of the far-field pattern with $n_{pan}=5N$, $N=1,2,3,\ldots, 24$, respectively.
}\label{fig1}
\end{figure}

FIG. \ref{fig1} shows the relative error
${|u^\infty(\hat{x};y)-u^\infty_{app}(\hat{x};y)|}/{|u^\infty(\hat{x};y)|}$
of the far-field pattern with $n_{pan}=5N$, $N=1,2,3,\ldots, 24$, respectively.
It can be seen in FIG. \ref{fig1} that the RCIP method is very stable and capable of producing
an accurate value of the far-field pattern if the parameter $n_{pan}$ is large enough.

\section{Uniqueness of the inverse problem}\label{sec3+}
\setcounter{equation}{0}

In this section, we establish uniqueness in the inverse scattering problem of reconstructing
the locally rough surface from the far-field pattern associated with the incident plane wave at
a fixed wave number $k>0$. Precisely, we will prove that the locally rough surface can be uniquely
determined by the far field patterns $u^\infty(\hat{x};d)$ for all $\hat{x}\in\Sp^1_+,\;d\in\Sp^1_-$
at a fixed wave number $k>0$. Our proof follows the idea in \cite{ColtonKress2013}.
We now state the following lemma which can be seen as an extension of Rellich's lemma
for the case of bounded obstacles (see, e.g., \cite[Lemma 2.12]{ColtonKress2013}).

\begin{lemma}\label{le1_nr}
Assume that $u^s\in H^1_{loc}(\R^2_+\ba\ov{B}_{R})$ satisfies the Helmholtz equation (\ref{eq1})
in $\R^2_+\ba\ov{B}_{R}$, the Sommerfeld radiation condition (\ref{rc}) and the Neumann boundary
condition ${\pa u^s}/{\pa\nu}=0$ on $\G\ba\ov{\G}_R$. Assume further that the far-field pattern
$u^\infty$ of the scattering solution $u^s$ vanishes on $\Sp^1_+$.
Then $u^s$ vanishes in $\R^2_+\ba\ov{B}_{R}$.
\end{lemma}

\begin{proof}
We extend $u^s$ to a function $\wid{u}^s$ in $\R^2\ba\ov{B}_R$ such that $\wid{u}^s(x)=u^s(x)$
in $\R^2_+\ba\ov{B}_{R}$, $\wid{u}^s(x)=u^s(x^{re})$ in $\R^2_-\ba\ov{B}_{R}$, where $x^{re}$ is
defined in Section \ref{sec2}. Then, by reflection $\wid{u^s}$ satisfies the Helmholtz
equation (\ref{eq1}) in $\R^2\ba\ov{B}_{R}$ and the Sommerfeld radiation condition (\ref{rc})
uniformly for all directions $\hat{x}=x/|x|$ with $x\in\R^2\ba\ov{B}_{R}$. Further,
the far-field pattern $\wid{u}^\infty$ of the scattering solution $\wid{u}^s$ vanishes on $\Sp^1$.
Then, by Theorem 2.13 in \cite{ColtonKress2013} we obtain that $\wid{u}^s$ vanishes in $\R^2\ba\ov{B}_R$,
which completes the proof.
\end{proof}

To prove the uniqueness result for the inverse scattering problem, we also need to consider
the case with the incident field being the point source $u^i(x;y):=\Phi_k(x,y)$ with the source
position $y\in D_+$. Then the reflected field with respect to the flat plane $x_2=0$ is
$u^r(x;y)=\Phi_k(x,y^{re})$. For arbitrarily fixed $y\in D_+$ let $u^s(x;y)$ satisfy the problem (NP)
with the boundary data $f(x)=f(x;y):=-{\pa u^i(x;y)}/{\pa\nu(x)}-{\pa u^r(x;y)}/{\pa\nu(x)}$, $x\in\G$.
Then, by Theorem \ref{thm1-wp} $u^s(\cdot;y)\in H^1_{loc}(D_+)$ exists and is the unique scattering
solution associated with the locally rough surface and the incident point source $\Phi_k(x,y)$
at the source position $y\in D_+$ since $f(x)=0$ on $\G\ba\G_R$ and thus is compactly supported on $\G_R$,
and so, $u(x;y)=u^i(x;y)+u^r(x;y)+u^s(x;y)$ is the total field.
Now, let $u^s(x;d)$ be the scattering solution generated by the incident plane wave
$u^i(x;d)=\exp({ikd\cdot x})$ with the incident direction $d\in\Sp^1_-$, and so,
$u(x;d)=u^i(x;d)+u^r(x;d)+u^s(x;d)$ is the total field.
Further, denote by $u^\infty(\hat{x};y)$ the far-field pattern of $u^s(x;y)$, where $\hat{x}\in S_+$.
Then we have the following mixed reciprocity relation.

\begin{lemma}\label{le2_nr}
$u^\infty(d;y)=\g u(y;-d),\;\forall y\in D_+,\;d\in\Sp^1_+$, where $\g=e^{i\pi/4}/\sqrt{8k\pi}$.
\end{lemma}

\begin{proof}
Without loss of generality, we may assume that $y\in D_{R}^+$. Then, for $d\in\Sp^1_+$ and $\vep>0$
small enough we derive by Green's second theorem that
\be\label{eq6_nr}
&&\int_{\pa B_{R}^+}\left[u(x;y)\frac{\pa u(x;-d)}{\pa\nu(x)}
    -u(x;-d)\frac{\pa u(x;y)}{\pa\nu(x)}\right]ds(x)\no\\
&&\qquad+\int_{|x-y|=\vep}\left[u(x;y)\frac{\pa u(x;-d)}{\pa\nu(x)}
    -u(x;-d)\frac{\pa u(x;y)}{\pa\nu(x)}\right]ds(x)\no\\
&&\;=\int_{\G_R}\left[u(x;y)\frac{\pa u(x;-d)}{\pa\nu(x)}-u(x;-d)\frac{\pa u(x;y)}{\pa\nu(x)}\right]ds(x)
    =0,
\en
where $\nu$ is the unit normal on $\pa B^+_R$ pointing outside of $D^+_R$, on $\G_R$ directing into $D_+$
and on $\{x\in\R^2:|x-y|=\vep\}$ directing into $\{x\in\R^2:|x-y|\le\vep\}$, and use has been made of
the Neumann boundary condition on $\G_R$ of $u(x;y)$ and $u(x;-d)$ to get the last equality.

Since $u^s(x;y)$ is analytic in the neighborhood of $\{x\in\R^2:|x-y|=\vep\}$, by the singularity at
$x=y$ of $\Phi_k(x,y)$ and its derivative (see Section 3.4 of \cite{ColtonKress2013}) we obtain that
\be\label{eq8_nr}
-u(y;-d)=\lim_{\vep\rightarrow 0}\int_{|x-y|=\vep}\left[u(x;y)\frac{\pa u(x;-d)}{\pa\nu(x)}
-u(x;-d)\frac{\pa u(x;y)}{\pa\nu(x)}\right]ds(x).
\en
On the other hand, we extend $u(x;y),\;u^s(x;-d)$ into $\R^2_-\ba\ov{B}_R$ by reflection,
which are denoted again by $u(x;y),\;u^s(x;-d)$, respectively, such that
$u(x;y)=u(x^{re};y),\;u^s(x;-d)=u^s(x^{re};-d)$ for $x\in\R^2_-\ba\ov{B}_R$.
Then, by the reflection principle $u(x;y),\;u^s(x;-d)$ satisfy the Helmholtz equation (\ref{eq1})
in $\R^2\ba\ov{B}_{R}$ and the Sommerfeld radiation condition (\ref{rc}) uniformly for all direction
$\hat{x}=x/|x|$ with $x\in\R^2\ba\ov{B}_R$. Thus, by using Green's second theorem
in $\R^2_-\ba\ov{B}_R$ we derive immediately that
\be\label{eq7_nr}
0&=&\int_{\pa B_{R}}\left[u(x;y)\frac{\pa u^s(x;-d)}{\pa\nu(x)}
 -u^s(x;-d)\frac{\pa u(x;y)}{\pa\nu(x)}\right]ds(x)\no\\
&=&2\int_{\pa B_{R}^+}\left[u(x;y)\frac{\pa u^s(x;-d)}{\pa\nu(x)}
 -u^s(x;-d)\frac{\pa u(x;y)}{\pa\nu(x)}\right]ds(x).
\en
And, by using the formula (3.87) in \cite{ColtonKress2013} we further have
\be\no
u^\infty(d;y)&=&\g\int_{\pa B_{R}}\left[u(x;y)\frac{\pa u^i(x;-d)}{\pa\nu(x)}
  -u^i(x;-d)\frac{\pa u(x;y)}{\pa\nu(x)}\right]ds(x)\\ \no
&=&\g\int_{\pa B_{R}^+}\left[u(x;y)\left[\frac{\pa u^i(x;-d)}{\pa\nu(x)}
  +\frac{\pa u^r(x;-d)}{\pa\nu(x)}\right]\right]ds(x)\\ \label{eq9_nr}
&&-\g\int_{\pa B_{R}^+}\left[\left[u^i(x;-d)+u^r(x;-d)\right]\frac{\pa u(x;y)}{\pa\nu(x)}\right]ds(x),
\en
where we have used the fact that $u^i(x^{e};-d)=u^r(x;-d)$ to get the second equality.
Letting $\vep\to0$ in (\ref{eq6_nr}), using (\ref{eq8_nr}), (\ref{eq7_nr}) and (\ref{eq9_nr})
and noting that $u(x;d)=u^i(x;d)+u^r(x;d)+u^s(x;d)$ give the required result.
This completes the proof.
\end{proof}

We are now ready to prove the following uniqueness theorem.

\begin{theorem}
Assume that $u^\infty_1(\hat{x};d)$ and $u^\infty_2(\hat{x};d)$ are the far-field patterns of the
scattering solutions for the locally rough surfaces $\G_1$ and $\G_2$, respectively, corresponding
to the same incident wave field $u^i(x;d)=\exp({ikd\cdot x})$.
If $u^\infty_1(\hat{x};d)=u^\infty_2(\hat{x};d)$ for all $\hat{x}\in\Sp^1_+$ and $d\in\Sp^1_-$
with a fixed wave number $k>0$ then $\G_1=\G_2$.
\end{theorem}

\begin{proof}
Denote by $D_{j+}$ and $D_{j-}$ the unbounded domains above and below the locally rough surface $\G_j$,
respectively, $j=1,2$. Suppose $D_{1+}\neq D_{2+}$. Without loss of generality, we may assume that
there exists $z\in\G_1\cap D_{2+}$ such that $\ov{B_h(z)}\subset D_{2+}$ for $h>0$ small enough, where
$B_h(z):=\{x\in\R^2:|x-z|<h\}$. Define $z_j=z+(h/j)\nu(z),\;j\in\N$,
where $\nu(z)$ is the unit normal at $z\in\G_1$ directed into $D_{1+}$.

Let $G$ be the unbounded connected part of $D_{1+}\cap D_{2+}$.
Denote by $\G_{j,p}$ the local perturbation part of $\G_j$, $j=1,2,$ and assume that
$\G_{j,p}\subset B_R$ for some large enough $R>0$, $j=1,2.$
Since $u^\infty_1(\hat{x};d)=u^\infty_2(\hat{x};d)$ for all $\hat{x}\in\Sp^1_+$ and $d\in\Sp^1_-$,
and by Lemma \ref{le1_nr} we have that $u^s_1(x;d)=u^s_2(x;d)$ for all $x\in\R^2_+\ba\ov{B}_R$ and
$d\in\Sp^1_-$. By the unique continuation principle it then follows that $u^s_1(x;d)=u^s_2(x;d)$
for all $x\in G$ and $d\in\Sp^1_-$, and so, $u_1(x;d)=u_2(x;d)$ for all $x\in G$ and $d\in\Sp^1_-$.
Further, by Lemma \ref{le2_nr} we have $u^\infty_1(d;y)=u^\infty_2(d;y)$ for all $d\in\Sp^1_+$ and $y\in G$.
By Lemma \ref{le1_nr} and the unique continuation principle again we obtain that $u^s_1(x;y)=u^s_2(x;y)$
for all $x,y\in G$ with $x\neq y$. This, together with the continuity of $u^s_i,i=1,2$, yields
that $u^s_1(x;y)=u^s_2(x;y)$ for all $x,y\in G$. In particular, $u^s_1(x;z_j)=u^s_2(x;z_j)$, $j\in\N$,
for all $x\in G$.

Since $z_j$, $j\in\N$, have a uniform distance from $\G_2$, then, by the well-posedness of the scattering
problem associated with the incident point sources $u^i(x;z_j)=\Phi_k(x,z_j)$, $j\in\N$, we have that
$\|\pa u^s_2(\cdot;z_j)/\pa\nu\|_{H^{1/2}(\G_1\cap B_h(z))}\le C$ with $C$ independent of $j$, $j=1,2,\ldots$.
On the other hand, by the boundary condition on $\G_1$ we see that
\ben
\lim\limits_{j\rightarrow\infty}\|\frac{\pa u^s_2(\cdot;z_j)}{\pa\nu}\|_{H^{1/2}(\G_1\cap B_h(z))}
&=&\lim\limits_{j\rightarrow\infty}\|\frac{\pa u^s_1(\cdot;z_j)}{\pa\nu}\|_{H^{1/2}(\G_1\cap B_h(z))}\\
&=&\lim\limits_{j\rightarrow\infty}\|\frac{\pa\Phi_k(\cdot,z_j)}{\pa\nu}\|_{H^{1/2}(\G_1\cap B_h(z))}=\infty.
\enn
This is a contradiction, and thus, $\G_1=\G_2$. The proof is then completed.
\end{proof}

\section{Numerical solution of the inverse problem}\label{se5}

In this section, we develop an inversion algorithm to solve the inverse problem numerically:
Given the far-field pattern $u^\infty$ of the scattering solution $u^s$ to the scattering
problem (\ref{eq1})-(\ref{rc}) (or the problem (NP)) associated with a finite number of incident
plane waves $u^i$ with multiple wave numbers $k$ (or multiple frequencies $\om$),
to determine the unknown locally rough surface $\G$ (or the unknown surface profile $h$
which describes the locally rough surface $\G$).
In particular, we propose a Newton-type reconstruction algorithm with multi-frequency data
for this inverse problem and conduct numerical experiments to demonstrate the effectiveness
of the reconstruction algorithm.

\subsection{Reconstruction algorithm with multi-frequency data}

We begin with introducing some notations. Given the incident plane wave $u^i_k(x;d)=e^{ikx\cdot d}$
with $d\in\Sp^1_-$, let $u^s_k$ be the solution of the scattering problem (\ref{eq1})-(\ref{rc})
with respect to the surface profile function $h$ characterizing the locally rough surface $\G$.
Assume that $u_k^\infty(\hat{x};d)$ is the far-field pattern of the scattering solution $u^s_k$.
Define the far-field operator $F_{d,k}:\;C^2(\R)\mapsto L^2(\Sp^1_+)$, which maps the function
$h$ to $u_k^\infty(\hat{x};d)$, by
\be\label{eq35}
\left(F_{d,k}[h]\right)(\hat{x})=u^\infty_k(\hat{x};d),\quad\hat{x}\in\Sp^1_+.
\en
Here, we use the subscript $k$ to indicate the dependence of $u^i_k,u^s_k,u_k^\infty$ and $F_{d,k}$
on the wave number $k$.

Our reconstruction algorithm consists in solving the nonlinear and ill-posed equation (\ref{eq35})
for the unknown function $h$ with using the Newton-type iterative method. To this end, we need to
investigate the Frechet differentiability of $F_{d,k}$ at $h$.
Let $\triangle h\in C^2_{0,R}(\R):=\{h\in C^2(\R):\textrm{supp}(h)\subset(-R,R)\}$ be a small
perturbation and let $\G_{\triangle h}:=\{(x_1,h(x)+\triangle h(x)):x_1\in\R\}$ denote the
corresponding locally rough surface characterized by $h(x)+\triangle h(x)$.
Then $F_{d,k}$ is called Frechet differentiable at $h$ if there exists a linear bounded operator $F'_{d,k}:C^2_{0,R}(\R)\rightarrow L^2(\Sp^1_+)$ such that
\ben
\left\|F_{d,k}[h+\triangle h]-F_{d,k}[h]-F'_{d,k}[h;\triangle h]\right\|_{L^2(\Sp^1_+)}
=o\left(\|\triangle h\|_{C^2(\R)}\right).
\enn

For the Frechet differentiability of $F_{d,k}$, we have the following theorem.

\begin{theorem}\label{thm2}
Let $u(x;d)=u^i(x;d)+u^r(x;d)+u^s(x;d)$, where $u^s$ solves the problem (NP)
with the boundary data $f=-(\pa u^i/\pa\nu+\pa u^r/\pa\nu)$. If $h\in C^2_{0,R}(\R)$,
then $F_{d,k}$ is Frechet differentiable at $h$ with the derivative given
by $F'_{d,k}[h;\triangle h]=(u')^\infty$ for $\triangle h\in C^2_{0,R}(\R)$.
Here, $(u')^\infty$ is the far-field pattern of $u'$ which solves the problem (NP)
with the boundary data $f=(d/ds)[(\nu_2\triangle h)(du/ds)]+k^2(\nu_2\triangle h)u$,
where $d/ds$ is the derivative with respect to the arc length and $\nu_2$ is the second component
of the unit normal vector $\nu$ on $\G$ directed into the unbounded domain $D_+$.
\end{theorem}

\begin{proof}
The proof is similar to that of \cite[Theorem 4.3]{BaoGaoLi2011},
where the Frechet derivative of $F_{d,k}$ is derived for the case of cavity problems
with a Neumann boundary condition by using an variational method.
\end{proof}

With the aid of the Frechet derivative of the far-field operator $F_{d,k}$, we are now ready to describe
the Newton-type iteration method for solving the inverse problem. Let $d_l\in\Sp^1_-,\;l=1,2,\cdots,n_d$,
be the incident directions and let $k>0$ be the wave number.
Assume that $h^{app}$ is an approximation function of $h$. For $d=d_l$ consider the following linearized
equations for (\ref{eq35}):
\be\label{eq36}
\big(F_{d_{l},k}[h^{app}]\big)(\hat{x})+\big(F'_{d_{l},k}[h^{app};\triangle h]\big)(\hat{x})
=u^\infty_k(\hat{x};d_{l}),\quad l=1,2,\ldots,n_d,
\en
where $\triangle h$ is the update function to be determined. The Newton-type method consists in iterating
the equations (\ref{eq36}) by using the Levenberg-Marquardt algorithm (see, e.g., \cite{Hohage1999}).

In the numerical experiments, we use the noised far-field pattern
$u^\infty_{\delta,k}(\hat{x_j};d_l),\;j=1,2,\ldots,n_f,\;l=1,2,\ldots,n_d$,
as the measurement data satisfying that
\ben
\left\|u^\infty_{\delta,k}(\cdot;d_l)-u^\infty_k(\cdot;d_l)\right\|_{L^2(\Sp^1_+)}
\le\delta\left\|u^\infty_k(\cdot;d_l)\right\|_{L^2(\Sp^1_+)},
\quad l=1,2,\ldots,n_d,
\enn
where $\delta>0$ is the noisy ratio and the observation directions $\hat{x_j},\;j=1,2,\ldots,n_f$,
are the equidistant points on $\Sp^1_+$.
In practical computation, $h^{app}$ has to be taken from a finite-dimensional subspace $R_M$,
where $R_M=\textrm{span}\{\phi_{1,M},\phi_{2,M},\cdots,\phi_{M,M}\}$ is a subspace of $C^2_{0,R}(\R)$
with $\phi_{j,M},\;j=1,2,\ldots,M$, the spline basis functions supported in $(-R,R)$ (see remark \ref{re5}).
Then, by the strategy in \cite{Hohage1999} we seek an update function
$\triangle h=\sum^M_{i=1}\triangle a_i\phi_{i,M}$ in $R_M$ which solves the minimization problem
\be\label{eq37}
&&\min_{\triangle a_i}\left\{\sum_{l=1}^{n_d}\sum_{j=1}^{n_f}\Big|\big(F_{d_{l},k}[h^{app}]\big)(\hat{x}_j)
+\big(F_{d_{l},k}'[h^{app};\triangle h]\big)(\hat{x}_j)-u^\infty_{\delta,k}(\hat{x_j};d_{l})\Big|^2\right.\no\\
&&\hspace{7cm}\left.+\beta\sum_{i=1}^M|\triangle a_i|^2\right\},
\en
where $\beta>0$ is chosen such that
\be\label{eq38}
\sum_{l=1}^{n_d}\sum_{j=1}^{n_f}\Big|\big(F_{d_{l},k}[h^{app}]\big)(\hat{x}_j)
+\big(F_{d_{l},k}'[h^{app};\triangle h]\big)(\hat{x}_j)-u^\infty_{\delta,k}(\hat{x_j};d_{l})\Big|^2\no\\
=\rho^2\sum_{l=1}^{n_d}\sum_{j=1}^{n_f}\Big|\big(F_{d_{l},k}[h^{app}]\big)(\hat{x}_j)
-u^\infty_{\delta,k}(\hat{x_j};d_{l})\Big|^2
\en
for a given constant $\rho\in(0,1)$. In computation, $\beta$ is determined by using
the bisection algorithm (see \cite{Hohage1999}). Note that $\beta$ is unique if certain
assumption is satisfied (see \cite{Hanke1997}). Then the approximation function
$h^{app}$ is updated by $h^{app}+\triangle h$. Further, define the error function
\ben
Err_k:=\frac{1}{n_d}\sum_{l=1}^{n_d}
\left[\sum\limits_{j=1}^{n_f}\left|\left(F_{d_{l},k}[h^{app}]\right)(\hat{x}_j)
-u^\infty_{\delta,k}(\hat{x_j};d_{l})\right|^2\right]^\half
\left[\sum\limits_{j=1}^{n_f}\left|u^\infty_{\delta,k}(\hat{x_j};d_{l})\right|^2\right]^{-1/2}.
\enn
Then the iteration is stopped if $Err_k<\tau\delta$, where $\tau>1$ is a fixed constant.

\begin{remark}\label{re5} {\rm
For a positive integer $M\in\N^+$ let $h=2R/(M+5)$ and $t_i=(i+2)h-R$.
Then the spline basis functions of $R_M$ can be defined as
$\phi_{i,M}(t)=\phi((t-t_i)/h),i=1,2,\ldots,M,$ where
\ben
\phi(t):=\sum^{k+1}_{j=0}\frac{(-1)^j}{k!}\left(\begin{array}{c}k+1\\
                                      j\end{array}\right)
       \left(t+\frac{k+1}{2}-j\right)^k_+
\enn
with $z^k_+=z^k$ for $z\geq0$ and $=0$ for $z<0$.
In this paper, we choose $k=4$, that is, $\phi$ is the quantic spline function.
Note that $\phi_{i,M}\in C^3(\R)$ with support in $(-R,R)$. See \cite{DeBoor} for details.
}
\end{remark}

\begin{remark}\label{re5+} {\rm
In the numerical experiments, the solution of the scattering problem will be computed numerically as follows.

1) The synthetic far-field data of the scattering problem is generated by using the RCIP method
introduced in Section \ref{se3-1} to solve the boundary integral equations (\ref{2.7}) and (\ref{eq14_nr}),
in which we choose the parameter $\rho=1$.

2) In each iteration, we use the boundary integral method proposed in Remark \ref{re2-9} to get
the numerical solution of the scattering problem to avoid the inverse crime and to reduce the
complexity of the computation. Precisely, assume that $u(x;d)=u^i(x;d)+u^r(x;d)+u^s(x;d),d\in\Sp^1_-$,
where $u^s$ is the solution of the problem (NP) with the boundary $\G$ replaced by $\G^{app}:=\{(x_1,x_2)\;|\;x_2=h^{app}(x_1),x_1\in\R\}$.
Then the scattered field $u^s$ is defined to be of the form (\ref{eq32}), where the density function
${\Phi}=(\varphi_1,\varphi_2)^T$ is the solution to the integral equation (\ref{eq31}).
Here, the boundary $\G$ appearing in (\ref{eq32}) and (\ref{eq31}) is also replaced by $\G^{app}$.
Similarly as in Section \ref{se3-1}, the integral equation (\ref{eq31}) can be numerically
solved to get the approximated values of $\Phi$ on the coarse mesh. Accordingly, the far-field pattern
$u^\infty$ of the scattered field $u^s$ can be computed with the formula (\ref{eq39}) (with $\G$ also
replaced by $\G^{app}$). Further, according to Theorem \ref{thm2} and the RCIP method introduced in
Section \ref{se3-1}, in order to compute the Frechet derivative at each iteration we need to compute
$(d/ds)[(\nu_2\triangle h)(du/ds)]+k^2(\nu_2\triangle h)u$ with $\triangle h\in R_M$
at the Gauss points of each panel for the coarse mesh of $\G^{app}_R:=\G^{app}\cap B_R$.
This can be done as follows. From the formula (\ref{eq32}) and the jump relation of the single-layer
potential, the scattered field $u^s$ on $\G^{app}_R$ is of the form
\be\label{eq41}
u^s(x)=(S_{k,\G^{app}_R\rightarrow\G^{app}_R}\varphi_1)(x)+(S_{k,\pa B^-_R\rightarrow\G^{app}_R}\varphi_2)(x),
\quad x\in\G^{app}_R.
\en
Then the approximated values $\rm\mathbf{u}^s_{coa}$ of the scattered field $u^s$ at the Gauss points for
the coarse mesh of $\G^{app}_R$ are computed as
\be\label{eq45}\rm
\mathbf{u}^s_{coa}=\mathbf{Q}\mathbf{L}_{fin}\mathbf{\Phi}_{fin},
\en
where $\mathbf{Q}$ is a restriction operator which performs panel-wise $15$-degree polynomial interpolation
in parameter from a grid on the fine mesh to a grid on the coarse mesh,
$\rm\mathbf{L}_{fin}\mathbf{\Phi}_{fin}$ is the discretization of the formula (\ref{eq41}) on the fine mesh,
and $\rm\mathbf{\Phi}_{fin}$ is the approximated values of $\mathbf{\Phi}$ on the fine mesh which can be
obtained by using a post-processor proposed in \cite[Section 7]{Helsing2009}.
Note that the approximated values $\rm\mathbf{u}^s_{coa}$ in (\ref{eq45}) can also
be obtained by using the method in \cite[Section 5]{Helsing2009}.
In this manner, the approximated values of the total field $u$ at the Gauss points for the coarse mesh of
$\G^{app}_R$ are thus obtained from (\ref{eq45}). Moreover, on each panel for the coarse mesh of $\G^{app}_R$,
we use the $15$-degree polynomial interpolation to approximate the total field $u$ on $\G^{app}_R$.
Finally, the first and second derivatives of the interpolation polynomial can be computed to obtain
the approximated values of $(d/ds)[(\nu_2\triangle h)({du}/{ds})]+k^2(\nu_2\triangle h)u$ on $\G^{app}_R$.
}
\end{remark}

Now we are ready to propose the Newton-type iteration algorithm for the inverse problem.
Motivated by \cite{BaoLin2011}, we use multi-frequency far-field data in order to get an accurate
reconstruction of the surface profile $h$. Generally speaking, the lower frequency data are used to
get the main profile of the boundary, whereas the higher frequency data can be applied to get a refined
reconstruction. The step of the algorithm is given as follows.

\begin{algorithm}\label{alg1}
Given the far-field data $u^\infty_{\delta,k_m}(\hat{x_j};d_{l}),j=1,2,\ldots,$
$n_f,l=1,2,\ldots,n_d,$ $m=1,2,\ldots,N$, with $\hat{x_j}\in\Sp^1_+,\;d_l\in\Sp^1_-$
and $k_1<k_2<\cdots<k_N$.
\begin{enumerate}
\item Let $h^{app}$ be the initial guess of $h$. Set $i=0$ and go to Step 2.
\item Set $i=i+1$. If $i>N$, then stop the iteration; otherwise, set $k=k_i$ and go to Step 3.
\item If $Err_k<\tau\delta$, go to Step 2; otherwise, go to Step 4.
\item Solve (\ref{eq37}) with the strategy (\ref{eq38}) to get an updated function $\triangle h$.
Let $h^{app}$ be updated by $h^{app}+\triangle h$ and go to Step 3.
\end{enumerate}
\end{algorithm}

\subsection{Numerical examples of the inverse problem}

In this section, several numerical experiments are carried out to demonstrate the effectiveness of
the inversion algorithm. We begin with the following assumptions in all numerical experiments.

\begin{enumerate}[1)]
\item
In each example we use multi-frequency data with the wave numbers $k=1,2,3,\ldots,N,$ where $N$ is the total
number of frequencies. Further, we choose $N\leq36$ and the auxiliary curve $\wid{\G}$
is taken to be a circle with radius $r=0.1$. It is easily seen that $kr$ with $k=1,2,\ldots,36$,
are not the zeros of the Bessel function $J_n$ for any integers $n$. Then, by Remark \ref{re6} it is valid
to use the integral equation (\ref{2.7}) to solve the Neumann problem (NP) for all examples.
\item
To generate the synthetic data and to compute the Frechet derivative in each iteration, for the coarse mesh,
we choose the number of panels on each smooth component $n_{pan}$ to be the nearest integer of $0.6*k+18$
for each wave number $k$. And for the fine mesh, we choose the number of subdivisions to be $n_{sub}=30$.
\item
We measure the half-aperture (that is, the measurement angle is between $0$ and $\pi$) far-field pattern
with $200$ equidistant points. The noisy data $u^\infty_{\delta,k}$ are obtained as
$u^\infty_{\delta,k}=u^\infty_k+\delta\zeta\|u^\infty_k\|_{L^2(\Sp^1_+)}/\|\zeta\|_{L^2(\Sp^1_+)}$,
where $\zeta$ is a random number with $\Rt(\zeta),\I(\zeta)$ belonging to the standard normal distribution
$N(0,1)$. For all cases, we choose the noisy ratio $\delta=5\%$.
\item
We set the parameters $\rho=0.8$ and $\tau=1.5$.
\item
In each figure, we use solid line '-' and dashed line '- -' to represent the actual curve and the
reconstructed curve, respectively.
\item
In all examples, for the shape of the local perturbation of the infinite plane we assume that
$\textrm{supp}(h)\in(-1,1)$; we further choose $R=1$ and use the smooth curves which are not in $R_M$.
\end{enumerate}

\begin{figure}[htbp]
\centering
\subfigure{\includegraphics[width=2.2in]{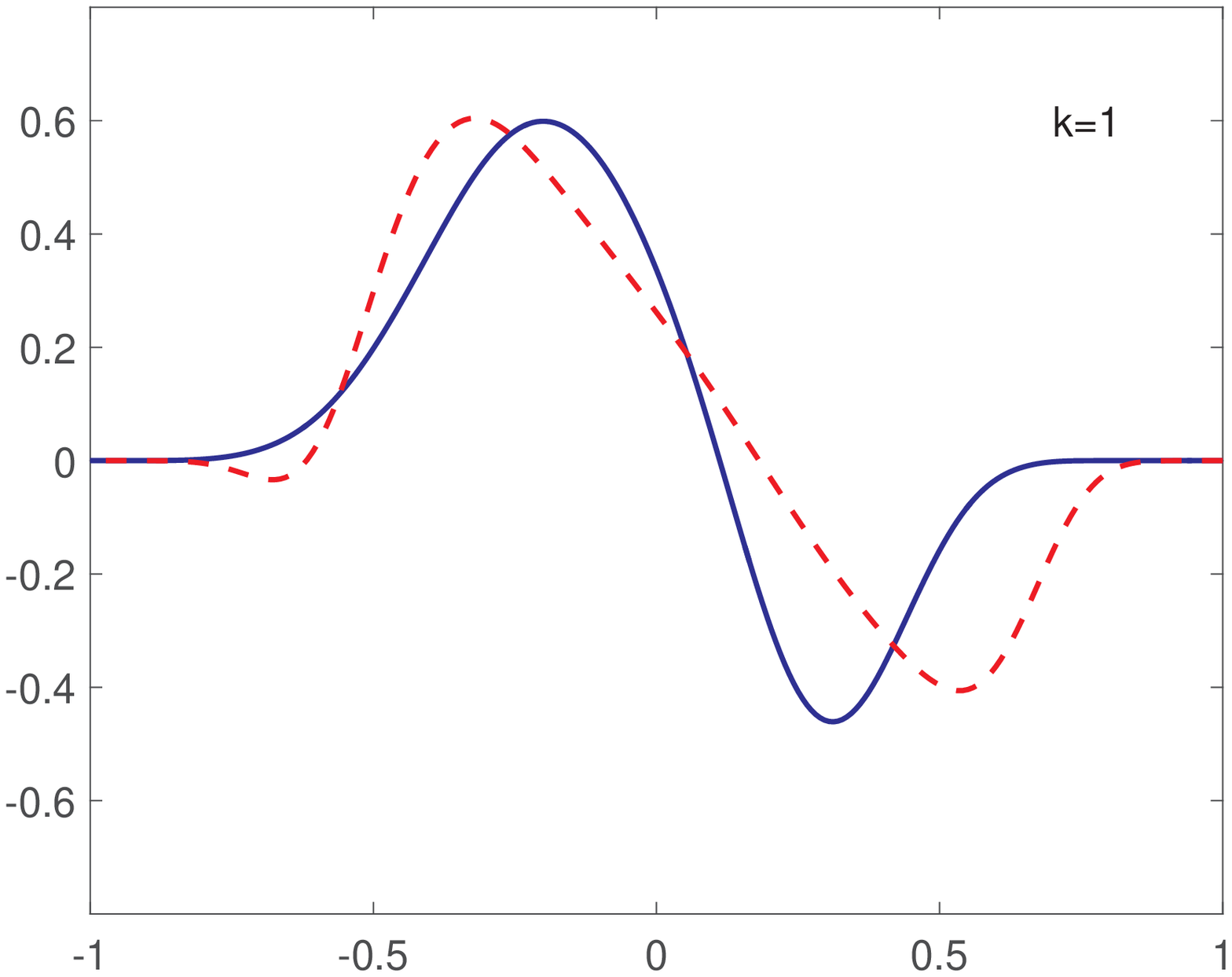}}
\hspace{0.2in}
\subfigure{\includegraphics[width=2.2in]{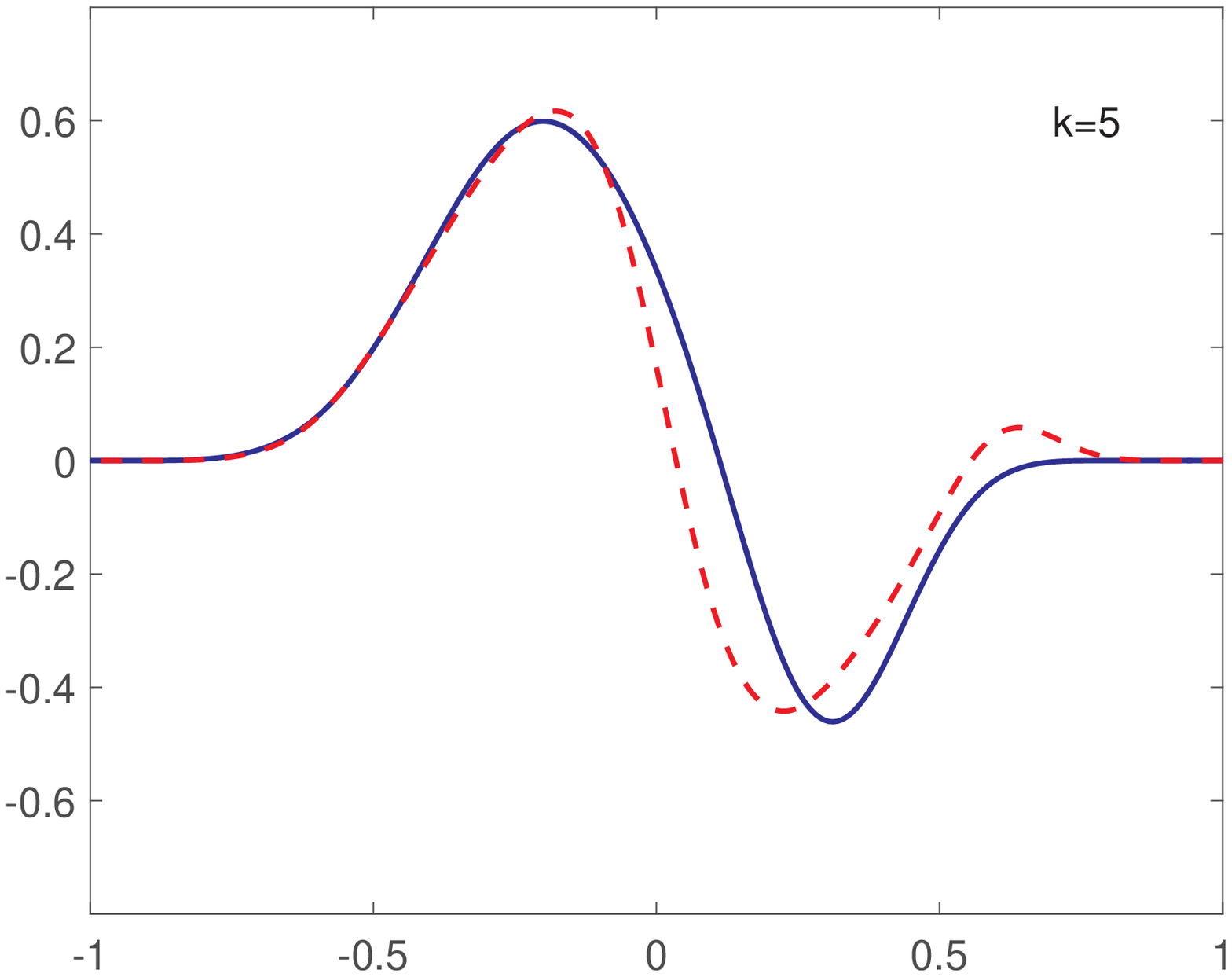}}
\subfigure{\includegraphics[width=2.2in]{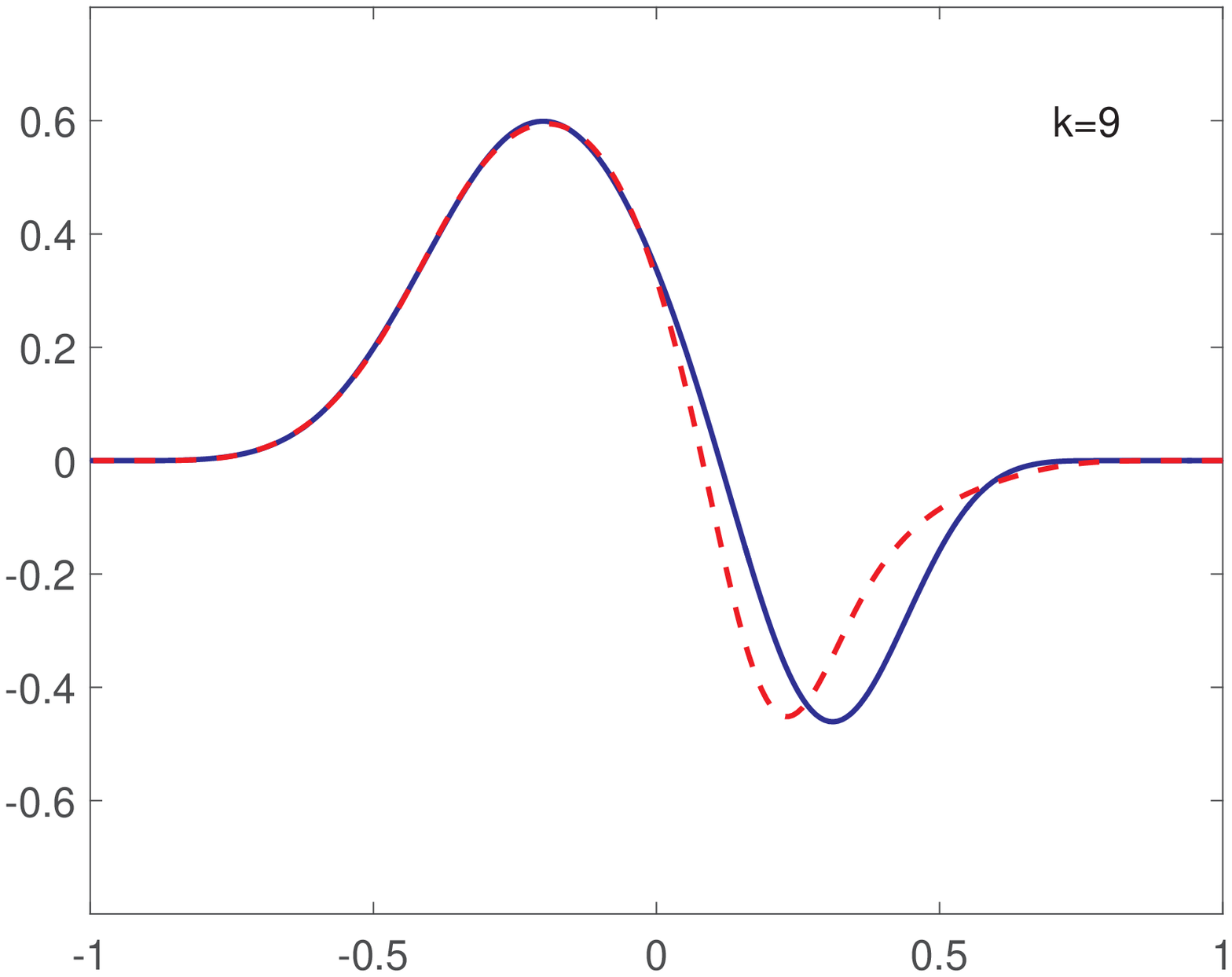}}
\hspace{0.2in}
\subfigure{\includegraphics[width=2.2in]{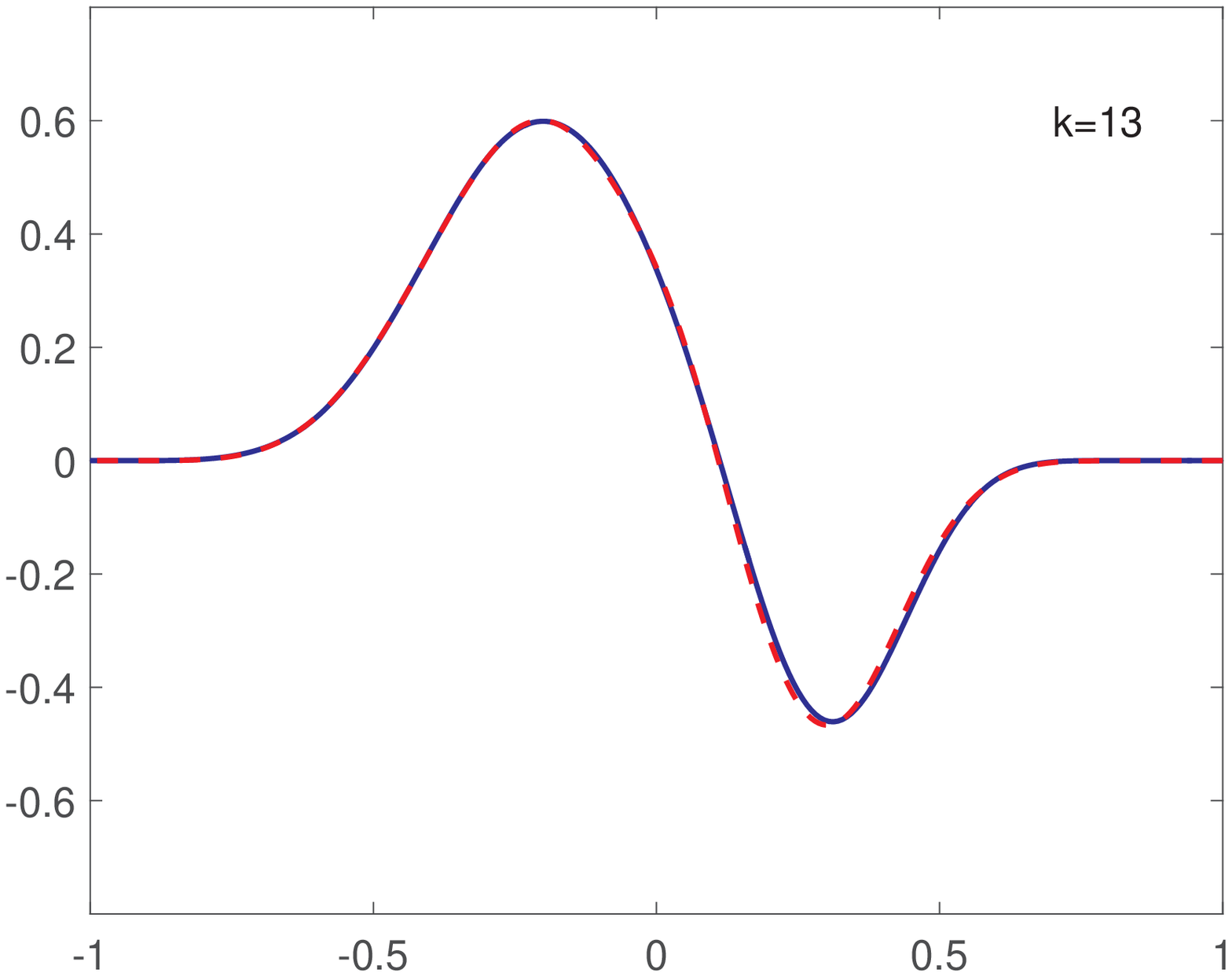}}
\caption{The reconstructed curve (dashed line) at $k=1,5,9,13,$ respectively, from $5\%$ noisy
far-field data for two incident fields $u^i(x;d_l),l=1,2$, with $d_1=(\cos(-\pi/3),\sin(-\pi/3))$
and $d_2=(\cos(-2\pi/3),\sin(-2\pi/3))$, where the real curve is denoted by the solid line.
}\label{fig2}
\end{figure}

\textbf{Example 1.} We first consider the surface profile
$h(x_1)=\phi(({x_1+0.2})/{0.3})-0.8\phi(({x_1-0.3})/{0.2}),$
where $\phi$ is defined in Remark \ref{re5}. We use the measured far-field data for two incident fields
$u^i(x;d_l),l=1,2$, with $d_1=(\cos(-\pi/3),\sin(-\pi/3))$ and $d_2=(\cos(-2\pi/3),\sin(-2\pi/3))$.
The auxiliary curve $\wid{\G}$ is chosen to be the circle centered at $[-0.3,-0.4]$ and of radius $r=0.1$.
For the inverse problem, we choose the number of the spline basis functions to be $M=40$ and the total number
of frequencies to be $N=13$. Figure \ref{fig2} shows the reconstructed curves at $k=1,5,9,13,$ respectively.

\begin{figure}[htbp]
\centering
\subfigure{\includegraphics[width=2.2in]{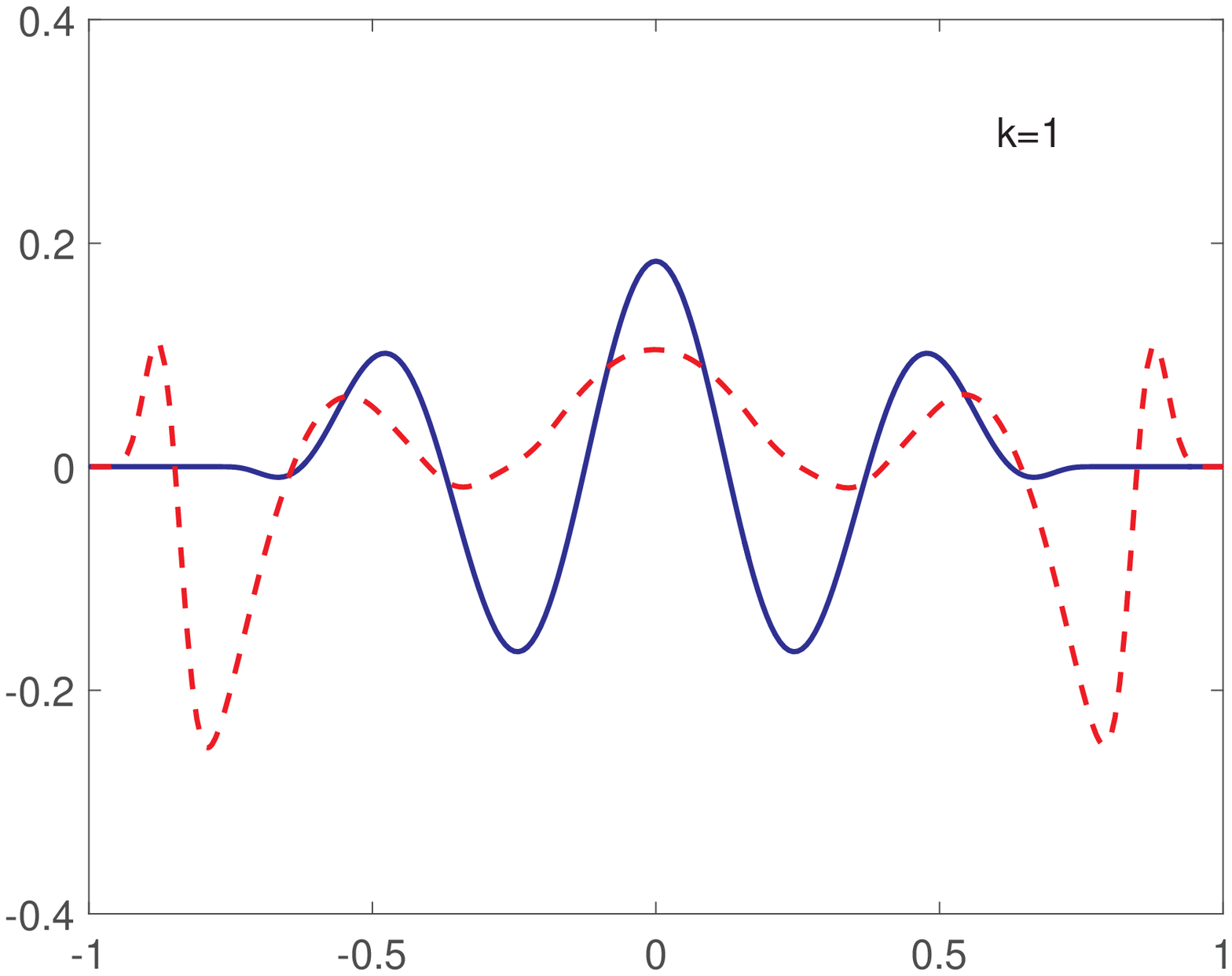}}
\hspace{0.2in}
\subfigure{\includegraphics[width=2.2in]{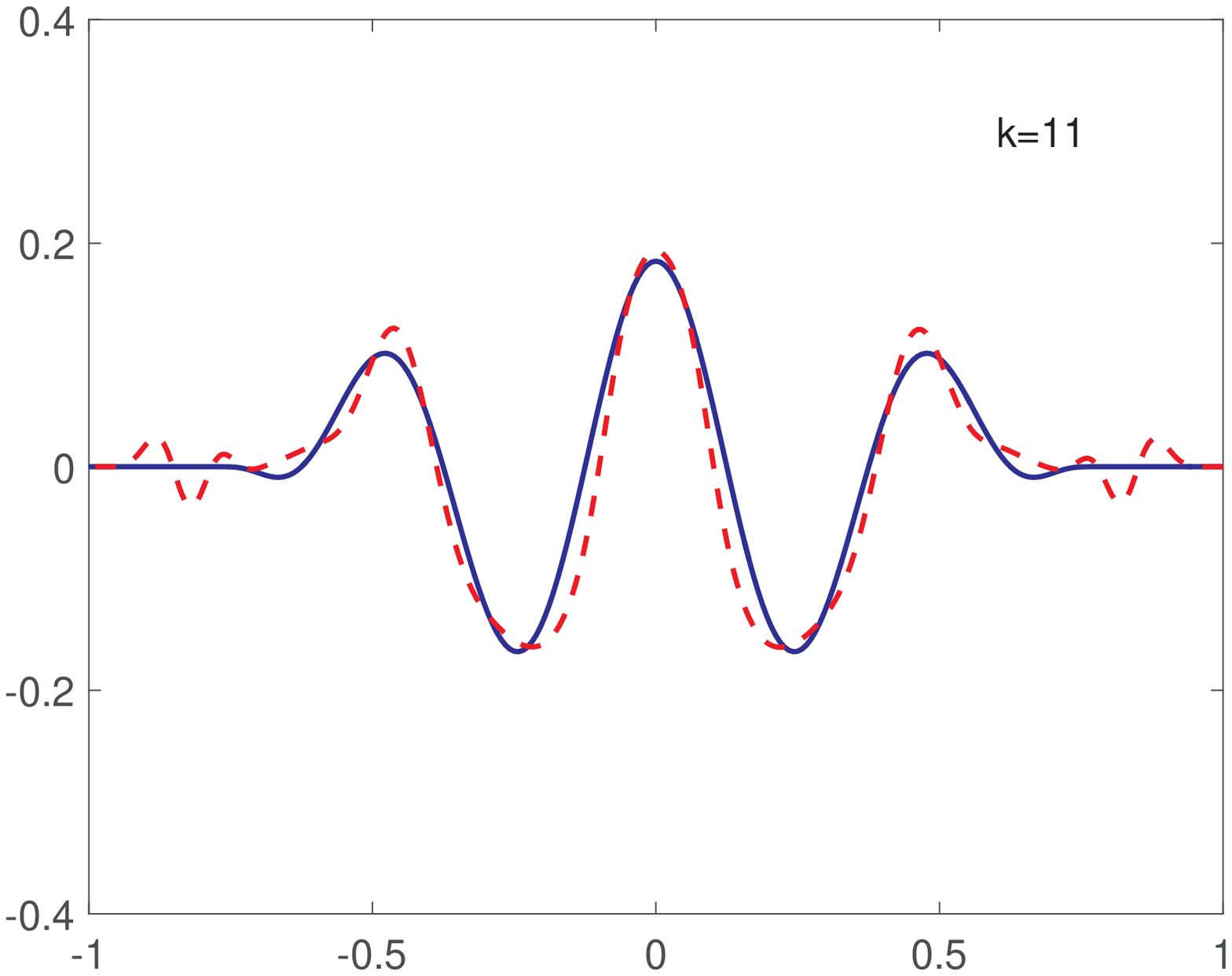}}
\subfigure{\includegraphics[width=2.2in]{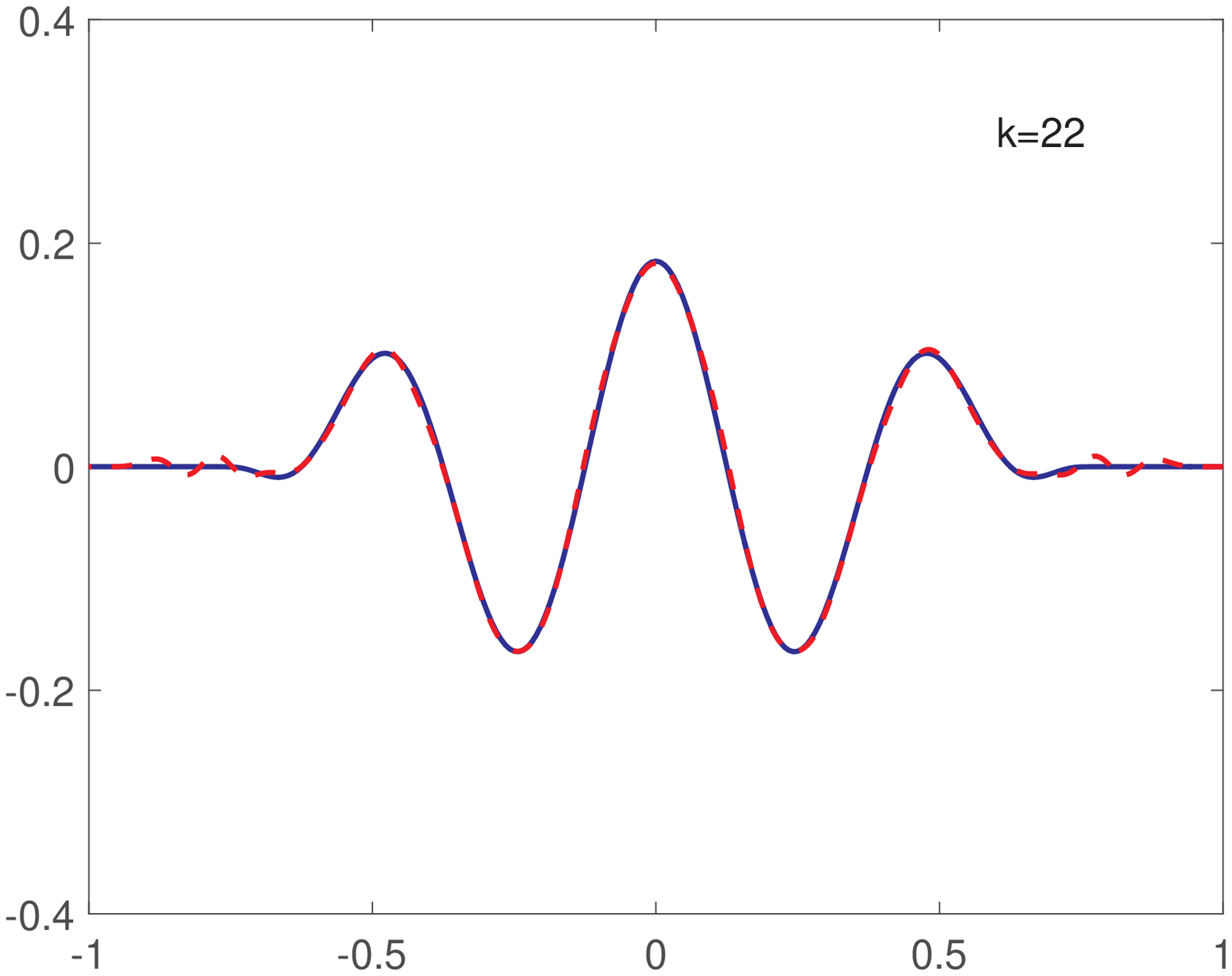}}
\hspace{0.2in}
\subfigure{\includegraphics[width=2.2in]{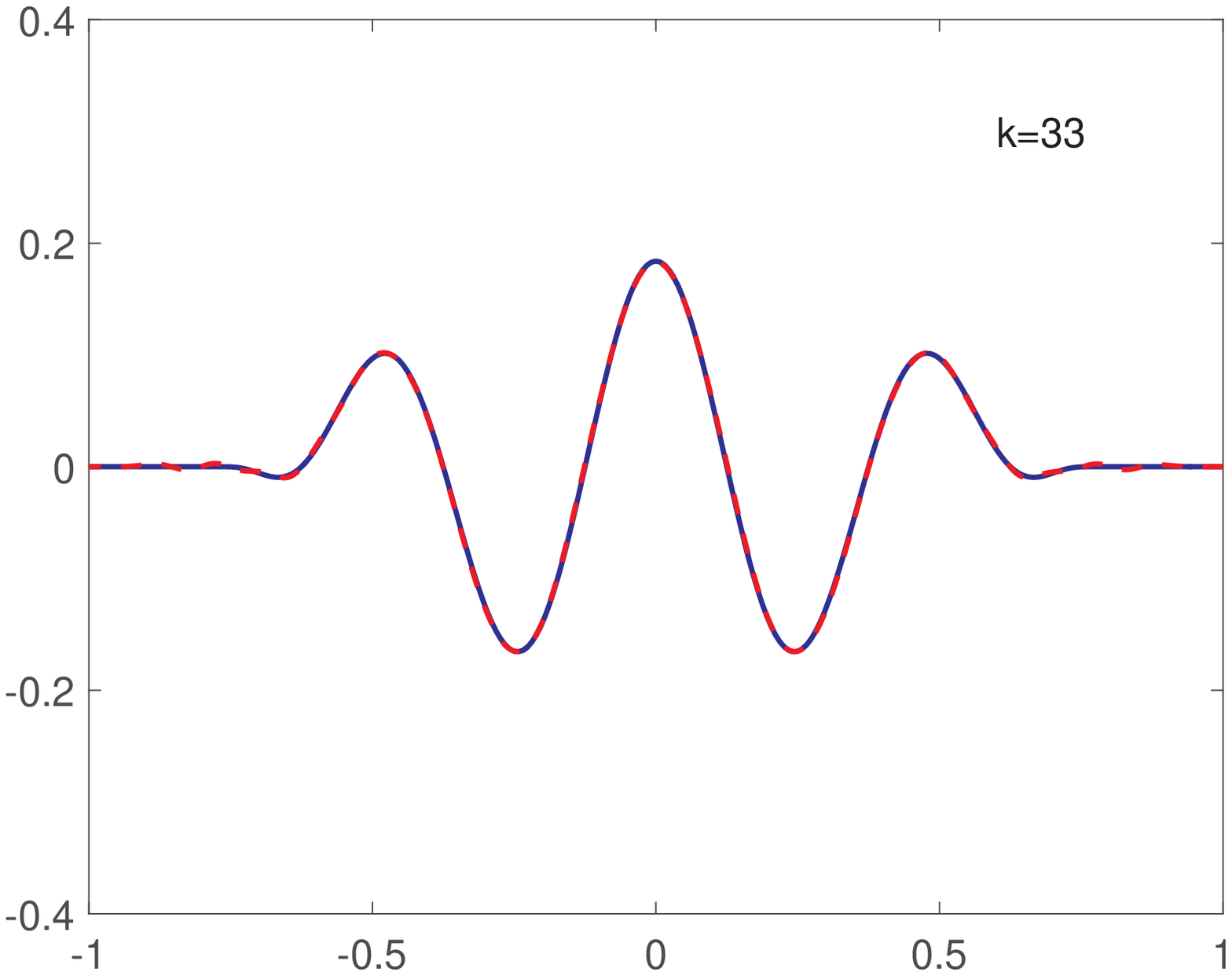}}
\caption{The reconstructed curve (dashed line) at $k=1,11,22,33,$ respectively, from $5\%$ noisy
far-field data for two incident fields $u^i(x;d_l),l=1,2$, with $d_1=(\cos(-\pi/3),\sin(-\pi/3))$
and $d_2=(\cos(-2\pi/3),\sin(-2\pi/3))$, where the real curve is denoted by the solid line.
}\label{fig3}
\end{figure}

\textbf{Example 2.} We now consider the surface profile defined by
\ben
h(x_1)=\left\{\begin{array}{ll}
\ds 0.5\exp\left[16/(25x_1^2-16)\right]\cos(4\pi x_1), &|x_1|<4/5,\\
\ds 0, &|x_1|\geq4/5.
\end{array}\right.
\enn
We use the measured far-field data for two incident fields $u^i(x;d_l),l=1,2$, with
$d_1=(\cos(-\pi/3),\sin(-\pi/3))$ and $d_2=(\cos(-2\pi/3),\sin(-2\pi/3))$.
The auxiliary curve $\wid{\G}$ is chosen to be the circle centered at $[0,-0.6]$ and with radius $r=0.1$.
For the inverse problem, the number of the spline basis functions is taken to be $M=40$
and the total number of frequencies is taken to be $N=33$.
Figure \ref{fig3} shows the reconstructed curves at $k=1,11,22,33,$ respectively.

\begin{figure}[htbp]
\centering
\subfigure{\includegraphics[width=2.2in]{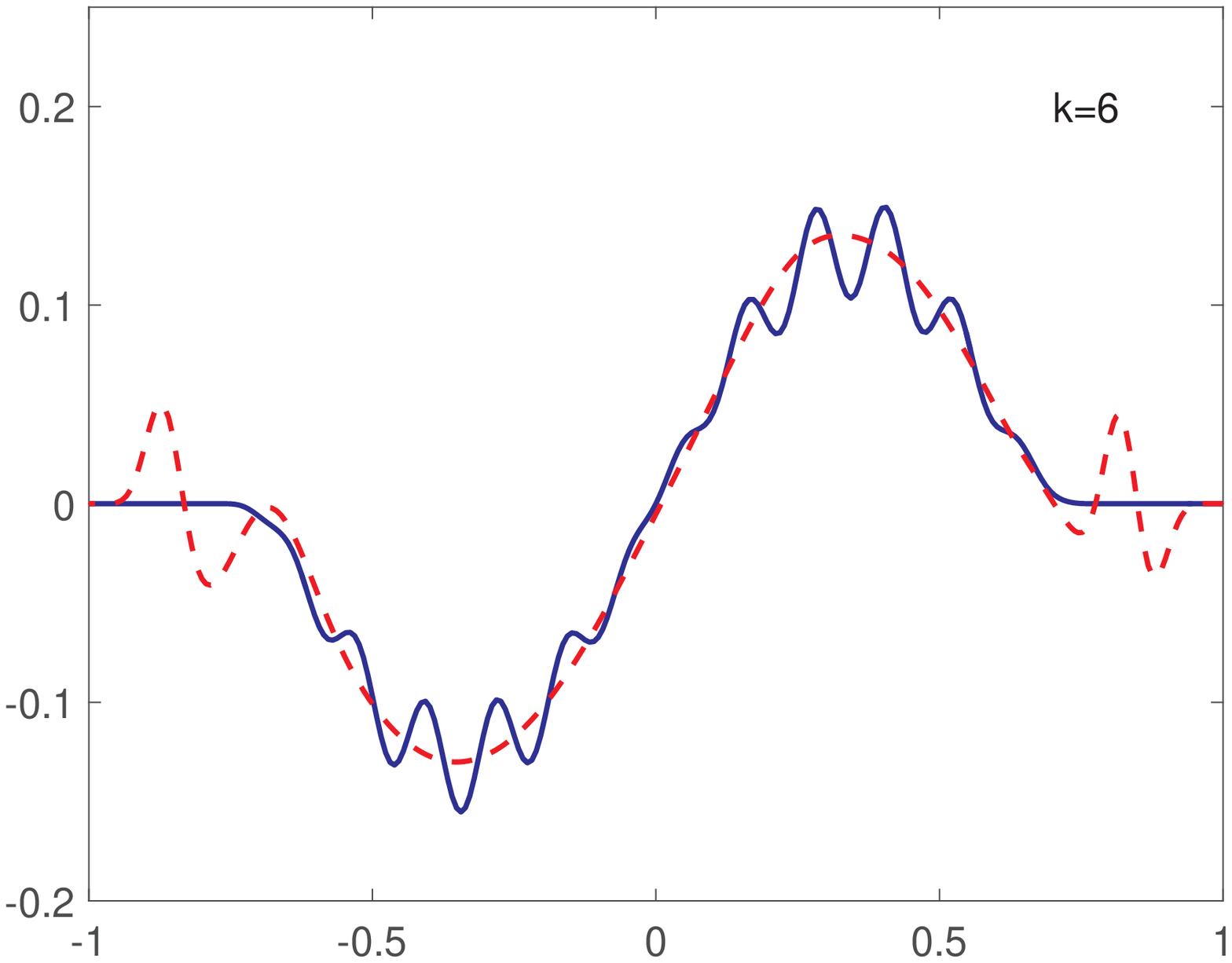}}
\hspace{0.2in}
\subfigure{\includegraphics[width=2.2in]{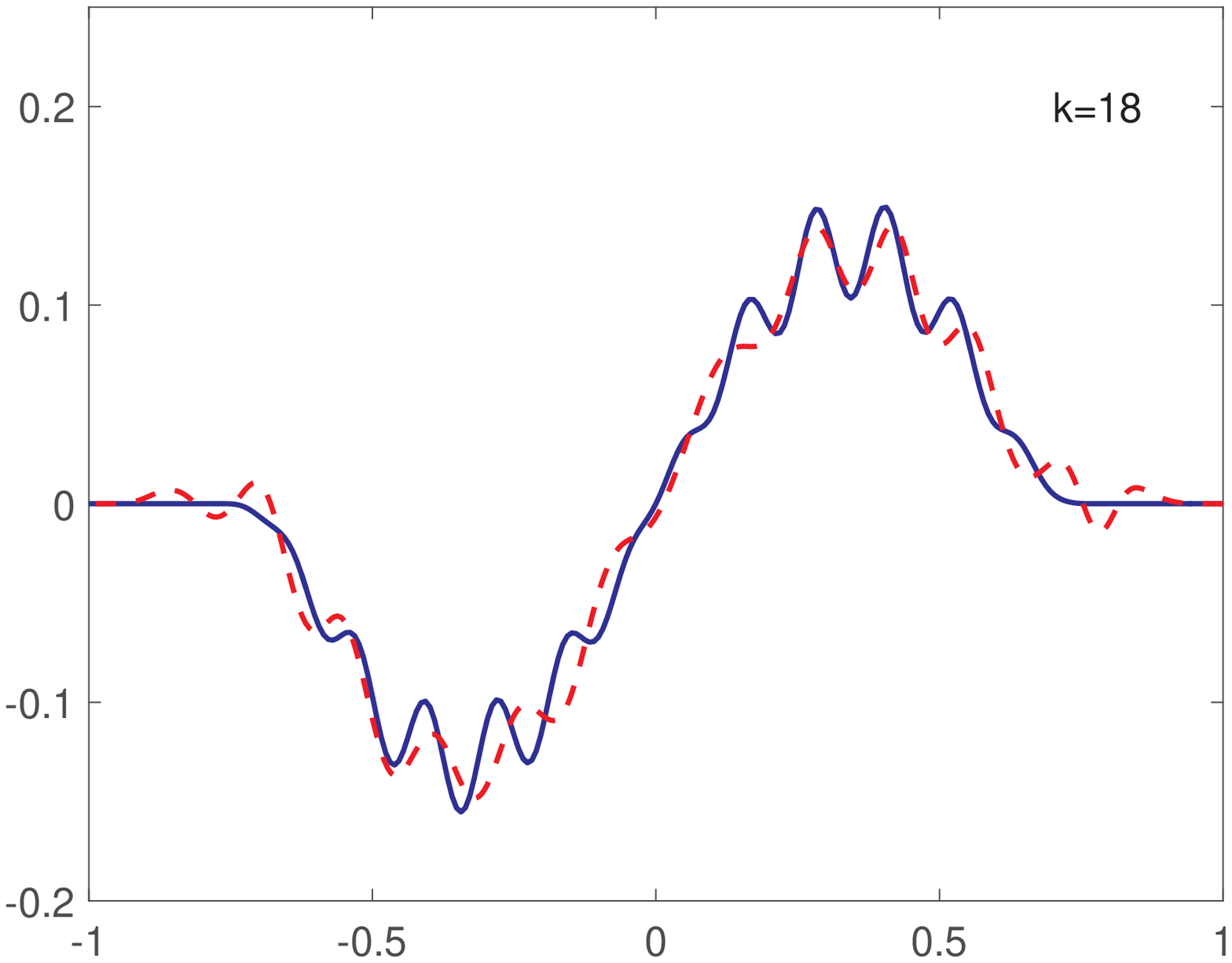}}
\subfigure{\includegraphics[width=2.2in]{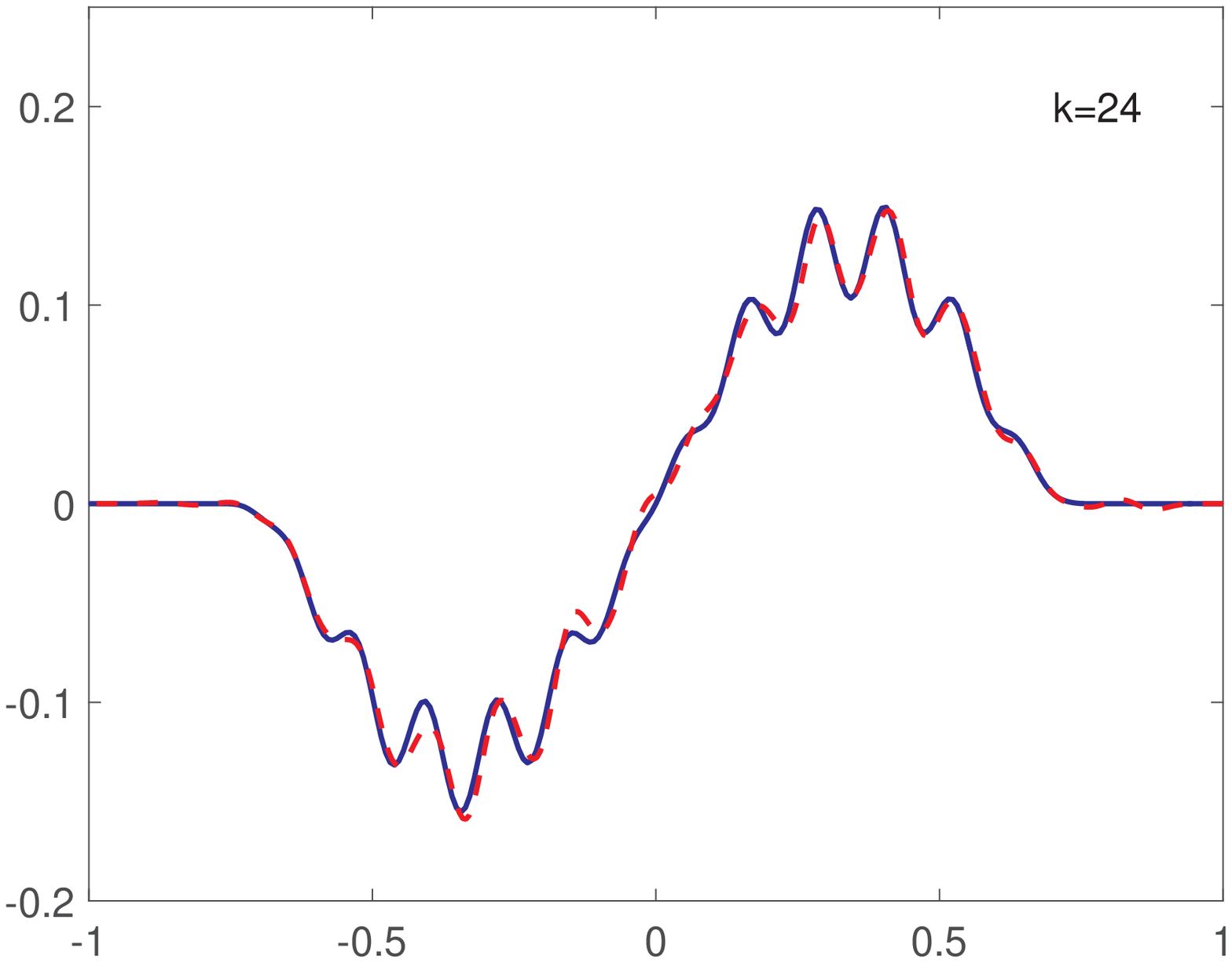}}
\hspace{0.2in}
\subfigure{\includegraphics[width=2.2in]{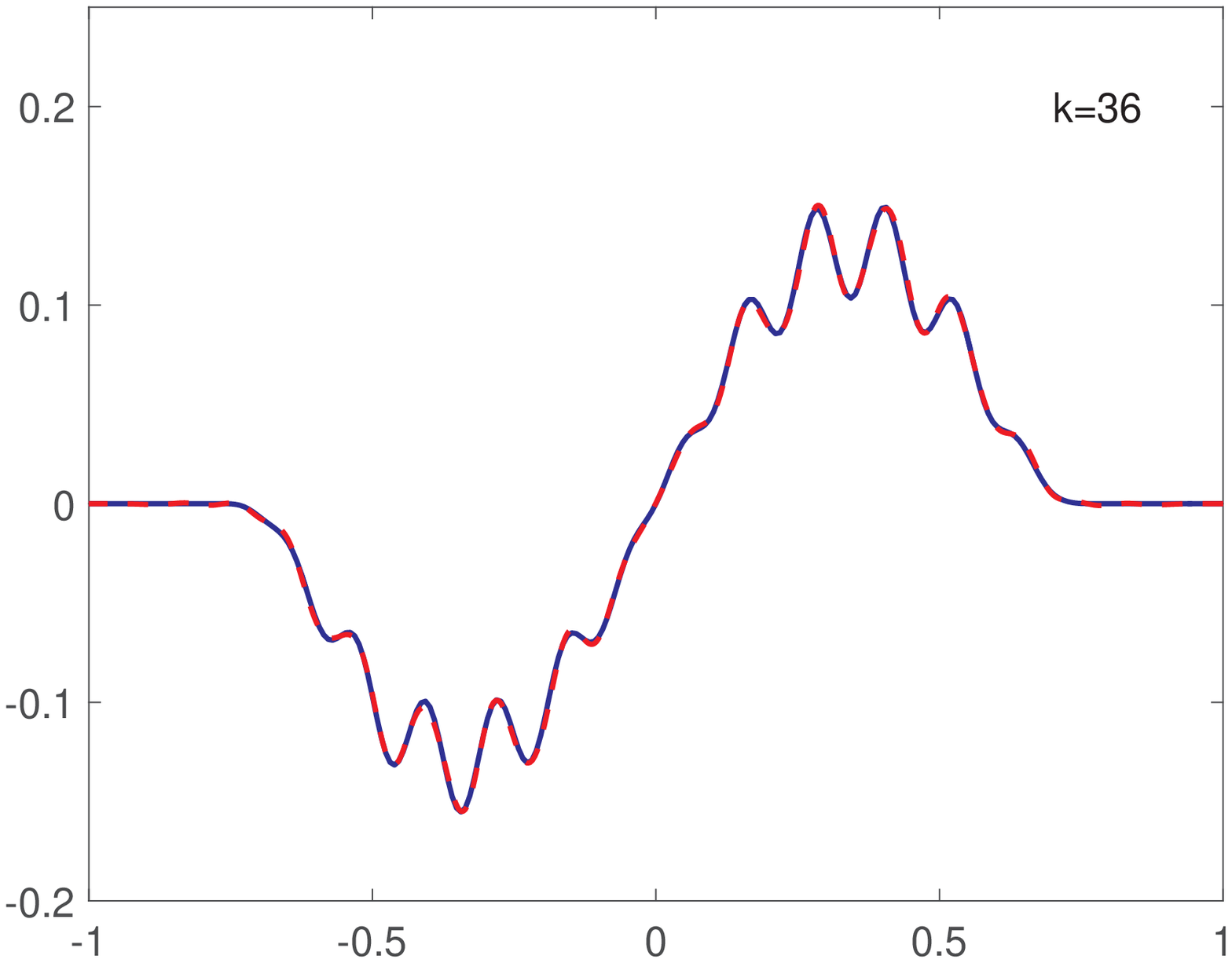}}
\caption{The reconstructed curve (dashed line) at $k=6,18,24,36,$ respectively, from $5\%$ noisy
far-field data for two incident fields $u^i(x;d_l),l=1,2$, with
$d_1=(\cos(-\pi/3),\sin(-\pi/3))$ and $d_2=(\cos(-2\pi/3),\sin(-2\pi/3))$,
where the real curve is represented by the solid line.
}\label{fig5}
\end{figure}

\textbf{Example 3 (multi-scale profile).} In this example, we consider a multi-scale profile given by
\ben
h(x_1)=\left\{\begin{array}{ll}
\ds \exp\left[16/(25x_1^2-16)\right]\left[0.5+0.1\sin(16\pi x_1)\right]\sin(\pi x_1), &|x_1|<4/5,\\
\ds 0, &|x_1|\geq4/5.
\end{array}\right.
\enn
This function is made up of two scales, that is, the macro-scale part
$0.5\exp\left[16/(25x_1^2-16)\right]\sin(\pi x_1)$ and the micro-scale part
$0.1\exp\left[16/(25x_1^2-16)\right]\sin(16\pi x_1)\sin(\pi x_1)$.
We use the measured far-field data for two incident fields $u^i(x;d_l),l=1,2$, with
$d_1=(\cos(-\pi/3),\sin(-\pi/3))$ and $d_2=(\cos(-2\pi/3),\sin(-2\pi/3))$.
We choose the auxiliary curve $\wid{\G}$ to be the circle with center at $[0,-0.4]$ and radius $r=0.1$.
For the inverse problem, the number of the spline basis functions is chosen to be $M=40$.
To capture the two-scale features of the profile, we choose the total number of frequencies to be $N=36$.
Figure \ref{fig5} presents the reconstructed curves at $k=6,18,24,36,$ respectively.
It is seen that the macro-scale feature is captured at the lower frequency $k=6$ (Figure \ref{fig5}, top left)
and the micro-scale feature is finally captured at the highest frequency $k=36$ (Figure \ref{fig5}, bottom right).

\section{Conclusion}\label{se6}

We considered the scattering problem by a locally rough surface on which the Neumann boundary condition
is imposed. An equivalent novel integral equation formulation was proposed for the direct scattering problem,
and a fast and efficient numerical algorithm was further given to solve the integral equation,
based on the RCIP Method previously introduced by J. Helsing.
The far-field pattern of the scattering problem can then be computed accurately by using the fast and efficient
integral equation solver.
For the inverse problem, we proved that the locally rough surface is uniquely determined from a knowledge
of the far-field pattern associated with incident plane waves. Further, a Newton-type iteration algorithm
with multi-frequency far-field data is also developed to reconstruct the locally rough surface,
where the proposed novel integral equation is applied to solve the direct scattering problem in each iteration.
Numerical experiments demonstrate that our inversion algorithm is stable and accurate even for the case of
multiple scale surface profiles. It is interesting to extend our method to the three-dimensional Helmholtz
and Maxwell equations which are more challenging.

\section*{Acknowledgements}

We thank Professor Johan Helsing for useful discussions on the RCIP method and especially
for his open Matlab codes which are very helpful. This work was supported by the NNSF of China
under grants 11871416, 11871466, and 91630309, and NSF of Shandong Province of China under Grant No. ZR2017MA044.

\end{document}